\documentclass[reqno,12pt,a4paper,final]{amsart}
\usepackage{amsmath,amsfonts,amssymb,amsthm,amstext,eucal}
\usepackage{bbold}
\usepackage{array}
\usepackage[latin1]{inputenc}
\newcommand{\half}{\tfrac{1}{2}}

\newcommand{\fd}{\mathfrak{d}}
\newcommand{\fg}{\mathfrak{g}}
\newcommand{\fh}{\mathfrak{h}}

\newcommand{\fS}{\mathfrak{S}}

\newcommand{\fso}{\mathfrak{so}}

\newcommand{\fsu}{\mathfrak{su}}
\newcommand{\fu}{\mathfrak{u}}

\newcommand{\SO}{\mathrm{SO}}
\renewcommand{\O}{\mathrm{O}}

\newcommand{\SU}{\mathrm{SU}}
\newcommand{\U}{\mathrm{U}}
\newcommand{\CC}{\mathbb{C}}
\newcommand{\EE}{\mathbb{E}}

\newcommand{\PP}{\mathbb{P}}
\newcommand{\RR}{\mathbb{R}}
\newcommand{\VV}{\mathbb{V}}

\DeclareMathOperator{\AdS}{AdS}

\DeclareMathOperator{\Gr}{Gr}
\DeclareMathOperator{\id}{id}

\newcommand{\repre}[1]{\boldsymbol{#1}}
\ifx\Sp\UnDeFiNeD
\newcommand{\Sp}{\mathrm{Sp}}
\else
\renewcommand{\Sp}{\mathrm{Sp}}
\fi
\newcommand{\MUNCH}[1]{\relax}
\allowdisplaybreaks[1]
\newtheorem{thm}{Theorem}
\newtheorem{conj}{Conjecture}
\begin{document}
\title{Plücker-type relations for orthogonal planes}
\author[Figueroa-O'Farrill]{José Figueroa-O'Farrill}
\address{School of Mathematics, University of Edinburgh, Scotland, UK}
\email{j.m.figueroa@ed.ac.uk}
\author[Papadopoulos]{George Papadopoulos}
\address{Department of Mathematics, King's College, London, England,
  UK}
\email{gpapas@mth.kcl.ac.uk}
\thanks{EMPG-02-17}
%\date{\today}
\begin{abstract}
  We explore a Plücker-type relation which occurs
  naturally in the study of maximally supersymmetric solutions of
  certain supergravity theories.  This relation generalises at the
  same time the classical Plücker relation and the Jacobi identity
  for a metric Lie algebra and coincides with the Jacobi identity of
  a metric $n$-Lie algebra.  In low dimension we present evidence for
  a geometric characterisation of the relation in terms of
  middle-dimensional orthogonal planes in euclidean or lorentzian
  inner product spaces.
\end{abstract}

\maketitle

\tableofcontents

\section{Introduction and main result}

The purpose of this note is to present a conjectural Plücker-style
formula for middle-dimensional orthogonal planes in real vector spaces
equipped with an inner product of euclidean or lorentzian signatures.
The formula is both a natural generalisation of the classical Plücker
formula and of the Jacobi identity for Lie algebras admitting an
invariant scalar product.  The formula occurs naturally in the study
of maximally supersymmetric solutions of ten-dimensional type IIB
supergravity and also in six-dimensional chiral supergravity.  We will
state the conjecture and then prove it for special cases which have
found applications in physics.  To place it in its proper mathematical
context we start by reviewing the classical Plücker relations.  For a
recent discussion see the paper \cite{EMPluecker} by Eastwood and
Michor.

\subsection{The classical Plücker relations}

The classical Plücker relations describe the projective embedding of
the grassmannian of planes.  Let $\VV$ be a $d$-dimensional vector
space (over $\RR$ or $\CC$, say) and let $\VV^*$ be the dual.  Let
$\Lambda^p \VV^*$ denote the space of $p$-forms on $\VV$ and
$\Lambda^p \VV$ the space of $p$-polyvectors on $\VV$.  We shall say
that a $p$-form $F$ is \emph{simple} (or \emph{decomposable}) if it
can be written as the wedge product of $p$ one-forms.  Every (nonzero)
simple $p$-form defines a $p$-plane $\Pi \subset \VV^*$, by declaring
$\Pi$ to be the span of the $p$ one-forms.  Conversely to such a
$p$-plane $\Pi$ one can associate a simple $p$-form by taking a basis
and wedging the elements together.  A different choice of basis merely
results in a nonzero multiple (the determinant of the change of basis)
of the simple $p$-form.  This means that the space of $p$-planes is
naturally identified with the subset of the projective space of the
space of $p$-forms corresponding to the rays of simple $p$-forms.  The
classical Plücker relations (see, e.g., \cite{GH,EMPluecker}) give
the explicit embedding in terms of the intersection of a number of
quadrics in $\Lambda^p \VV^*$.  Explicitly one has the following:

\begin{thm}
  A $p$-form $F \in \Lambda^p\VV^*$ is \emph{simple} if and only if
  for every ($p{-}1$)-polyvector $\Xi \in \Lambda^{p-1}\VV$,
  \begin{equation*}
    \iota_{\Xi} F \wedge F = 0~,
  \end{equation*}
  where $\iota_\Xi F$ denotes the one-form obtained by contracting
  $F$ with $\Xi$.
\end{thm}

Being homogeneous, these equations are well-defined in the projective
space $\PP\Lambda^p \VV^* \cong \PP^{\binom{d}{p}-1}$, and hence
define an algebraic embedding there of the grassmannian $\Gr(p,d)$ of
$p$-planes in $d$ dimensions.

The Plücker relations arise naturally in the study of maximally
supersymmetric solutions of eleven-dimensional supergravity
\cite{Bonn,FOPMax}.  Indeed, the Plücker relations for the $4$-form
$F_4$ in eleven-dimensional supergravity arise from the zero
curvature condition for the supercovariant derivative.  A similar
analysis for ten-dimensional type IIB supergravity \cite{FOPMax}
yields new (at least to us) Plücker-type relations, to which we now
turn.

\subsection{Orthogonal Plücker-type relations}

Let $\VV$ be a real vector space of finite dimension equipped with a
euclidean or lorentzian inner product $\left<-,-\right>$.  Let $F \in
\Lambda^p \VV^*$ be a $p$-form and let $\Xi \in \Lambda^{p-2}\VV$ be
an ($p{-}2$)-polyvector.  The contraction $\iota_\Xi F$ of $F$ with
$\Xi$ is a $2$-form on $\VV$ and hence gives rise to an element of the
Lie algebra $\fso(\VV)$.  If $\omega \in \Lambda^2\VV^* \cong
\fso(\VV)$, we will denote its action on a form $\Omega \in
\Lambda\VV^*$ by $[\omega, \Omega]$.  Explicitly, if $\omega = \alpha
\wedge \beta$, for $\alpha,\beta\in\VV^*$, then
\begin{equation*}
  [\alpha \wedge \beta, \Omega] = \alpha \wedge \iota_{\beta^\sharp}
  \Omega - \beta \wedge \iota_{\alpha^\sharp} \Omega~,
\end{equation*}
where $\alpha^\sharp\in\VV$ is the dual vector to $\alpha$ defined
using the inner product.  We then extend linearly to any $2$-form
$\omega$.

Let $F_1$ and $F_2$ be two simple forms in $\Lambda^p\VV^*$.  For the
purposes of this note we will say that $F_1$ and $F_2$ are
\textbf{orthogonal} if the $d$-planes $\Pi_i \subset \VV$ that they
define are orthogonal; that is, $\left<X_1,X_2\right>=0$ for all $X_i
\in \Pi_i$.  Note that if the inner product in $\VV$ is of lorentzian
signature then orthogonality does not imply that $\Pi_1 \cap \Pi_2 =
0$, as they could have a null direction in common.  If this is the
case, $F_i = \alpha \wedge \Theta_i$, where $\alpha$ is a null form
and $\Theta_i$ are orthogonal simple forms in a euclidean space in two
dimensions less.  Far from being a pathology, the case of null forms
plays an important role in the results of \cite{FOPMax} and is
responsible for the existence of a maximally supersymmetric plane wave
in IIB supergravity \cite{NewIIB}.

We now can state the following:

\begin{conj}
%\begin{itemize}

(i)~ Let $p\geq 2$ and $F\in\Lambda^p\VV^*$ be a $p$-form on an
  $d$-dimensional euclidean or lorentzian inner product space $\VV$,
  where $d=2p$ or $d=2p+1$.  For all ($p{-}2$)-polyvectors $\Xi \in
  \Lambda^{p-2}\VV$, the
  equation
  \begin{equation}
    \label{eq:conj}
    \left[ \iota_\Xi F, F \right] = 0
  \end{equation}
  is satisfied if and only if $F$ can be written as a sum of
  two orthogonal simple forms; that is,
  \begin{equation*}
    F= F_1 + F_2
  \end{equation*}
  where $F_1$ and $F_2$ are simple and $F_1 \perp F_2$.

  (ii)~ Let $p\geq 2$ and $F\in\Lambda^p\VV^*$ be a $p$-form on the
  euclidean or lorentzian vector space $\VV$ with dimension $p\leq d<2p$.
  The equation  \eqref{eq:conj} holds if and only if  $F$ is simple.
\end{conj}

Again the equation is homogeneous, hence its zero locus is
well-defined in the projective space of $\PP\Lambda^p \VV^* \cong
\PP^{\binom{d}{p}-1}$.

Relative to a basis $\{e_i\}$ for $\VV$ relative to which the inner
product has matrix $g_{ij}$, we can rewrite equation \eqref{eq:conj}
as
\begin{equation*}
  \sum_{k,\ell=1}^d g^{k\ell} F_{ki_1i_2\cdots i_{p-2}[j_1}
  F_{j_2j_3\cdots j_p]\ell} = 0~,
\end{equation*}
which shows that the ``if'' part of the conjecture follows trivially:
simply complete to a pseudo-orthonormal basis for $\VV$ the bases for
the planes $\Pi_i$, express this equation relative to that basis and
observe that every term vanishes.

Finally let us remark as a trivial check that both the equation
\eqref{eq:conj} and the conclusion of the conjecture are invariant
under the orthogonal group $\O(\VV)$.  A knowledge of the orbit
decomposition of the space of $p$-forms in $\VV$ under $\O(\VV)$ might
provide some further insight into this problem.

To this date the first part of the conjecture has been verified for the
following cases
\begin{itemize}
\item ($p{\leq}2$) both for euclidean and lorentzian signatures,
\item ($d{=}6$, $p{=}3$) both for euclidean and lorentzian signatures,
\item ($d{=}7$, $p{=}3$) for euclidean signature,
\item ($d{=}8$, $p{=}4$) for euclidean signature, and
\item ($d{=}10$, $p{=}$5) for euclidean and lorentzian signatures.
\end{itemize}
It is the latter case which is required in the investigation of
maximally supersymmetric solutions of ten-dimensional type IIB
supergravity \cite{FOPMax}, whereas the second case enters in the
case of six-dimensional $(1,0)$ supergravity \cite{CFOSchiral}.
The fourth is expected to have applications in eight-dimensional
supergravity theories.

The second part of the conjecture has been verified in the cases:
\begin{itemize}
\item ($p{\leq}2$) both for euclidean and lorentzian signatures,
\item ($d{<}6$, $p{=}3$) both for euclidean and lorentzian signatures,
\item ($d{<}8$, $p{=}4$) for euclidean signature.
\end{itemize}

There are two conditions in the hypothesis which seem artificial at
first:
\begin{itemize}
\item the restriction on the signature of the inner product, and
\item the restriction on the dimension of the vector space.
\end{itemize}
These conditions arise from explicit counterexamples for low $p$,
which we now discuss together with a Lie algebraic re-interpretation
of the identity \eqref{eq:conj}.

Before we proceed to explain these, let us remark that it might just
be the case that the restriction on the dimension of the vector space
is an artifact of low $p$.  We have no direct evidence of this, except
for the following.  We depart from the observation that the ratio of
the number of relations to the number of components of a $p$-form in
$d$ dimensions is $\binom{d}{p-2}$.  For fixed $p$ and large $d$, this
ratio behaves as $d^{p-2}$.  So for $p=2$, the ratio is one and for
$p=3$ grows linearly as $d$.  It is the latter case where the
counterexamples that justify the restriction on the dimensions will be
found.  For $p>3$ this ratio grows much faster and it is perhaps not
unreasonable to expect that the only solutions are those which verify
the conjecture.

\subsection{The case $p{=}2$}

Let us observe that for $p{=}2$ there are no equations, since
$[F,F]=0$ trivially in $\fso(\VV)$.  The conjecture would say that any
$F \in \fso(\VV)$ can be ``skew-diagonalised''.  In euclidean
signature this is true: it is the conjugacy theorem for Cartan
subalgebras of $\fso(\VV)\cong\fso(d)$.  The result also holds in
lorentzian signature; although it is more complicated, since depending
on the type of element (elliptic, parabolic or hyperbolic) of
$\fso(\VV) \cong \fso(1,d-1)$, it conjugates to one of a set of normal
forms, all of which satisfy the conjecture.

The conjecture does not hold in signature $(2,d)$ for any $d\geq 2$,
as a quick glance at the normal forms of elements of $\fso(2,d)$ under
$\O(2,d)$ shows that there are irreducible blocks of dimension higher
than $2$.  In other words, there are elements $\omega \in \fso(2,d)$
for which there is no decomposition of $\RR^{2,d}$ into 2-planes
stabilised by $\omega$.  A similar situation holds in signature
$(p,q)$ for $p,q > 2$, as can be gleaned from the normal forms
tabulated in \cite{Boubel}.

This justifies restricting the signature of the scalar product on
$\VV$ in the hypothesis to the conjecture.  The restriction on the
dimension of $\VV$ arises by studying the case $p{=}3$, to which we
now turn.

\subsection{The case $p{=}3$}

Let $F \in \Lambda^3 \VV^*$.  Using the scalar product $F$ defines a
linear map $[-,-]: \Lambda^2 \VV \to \VV$ by
\begin{equation}
  \label{eq:invmet}
  F(X,Y,Z) = \left<[X,Y],Z\right>~,\quad\text{for all $X,Y,Z\in\VV$.}
\end{equation}
The Plücker formula \eqref{eq:conj} in this case is nothing but the
statement that for all $X\in\VV$, the map $Y\mapsto [X,Y]$ should be a
derivation over $[-,-]$:
\begin{equation}
  \label{eq:jacobi}
  [X,[Y,Z]] = [[X,Y],Z] + [Y,[X,Z]]~.
\end{equation}
In other words, it is the Jacobi identity for $[-,-]$, turning $\VV$
into a Lie algebra, as the notation already suggests.  More is true,
however, and because of the fact that $F\in\Lambda^3\VV^*$, the metric
is invariant:
\begin{equation*}
  \left<[X,Y],Z\right> =   \left<X,[Y,Z]\right>~.
\end{equation*}
In other words, solutions of \eqref{eq:conj} for $p{=}3$ are in
one-to-one correspondence with Lie algebras admitting an invariant
nondegenerate scalar product.

We will show below (in two different ways) that the conjecture works
for $d\leq 7$, but the simple Lie algebra $\fsu(3)$ with the Killing
form provides a counterexample to the conjecture for $d=8$ (and also
for any $d>8$ by adding to it an abelian factor).  To see this,
suppose that the $3$-form $F$ associated to $\fsu(3)$ decomposed into
a sum\footnote{For general $p$ and $d$, there is no reason why $F$
  should break up as $F = F_1 + F_2$; in the general case we would
  have $F=F_1 + F_2 + \cdots$, where the $F_i$ are simple and mutually
  orthogonal.} $F=F_1 + F_2$ of orthogonal simple forms.  Each $F_i$
defines a $3$-plane in $\fsu(3)$.  Let $Z \in \fsu(3)$ be orthogonal
to both of these planes: such $Z$ exists because $\dim\fsu(3)=8$.
Then $\iota_Z F = 0$, and this would mean that for all $X,Y$,
$F(Z,X,Y) = \left<[Z,X],Y\right>=0$, so that $Z$ is central, which is
a contradiction because $\fsu(3)$ is simple.

\subsection{Metric $n$-Lie algebras}

There is another interpretation of the Plücker relation
\eqref{eq:conj} in terms of a generalisation of the notion of Lie
algebra.\footnote{We are grateful to Dmitriy Rumynin for making us
  aware of the existence of this concept.}

Let $p=n+1$ and $F \in \Lambda^{n+1} \VV^*$ and as we did for $p=3$
let us define a map $[\cdots]: \Lambda^n \VV \to \VV$ by
\begin{equation}
  \label{eq:n-Lie}
  F(X_1,X_2,\dots,X_{n+1}) = \left<[X_1,\dots,X_n],X_{n+1}\right>~.
\end{equation}
The relation \eqref{eq:conj} now says that for all
$X_1,\dots,X_{n-1}\in\VV$, the endomorphism of $\VV$ defined by $Y
\mapsto [X_1,\dots,X_{n-1},Y]$ is a derivation over $[\cdots]$; that
is,
\begin{equation}
  \label{eq:n-Jacobi}
  [X_1,\dots,X_{n-1},[Y_1,\cdots,Y_n]] = \sum_{i=1}^n
  [Y_1,\dots,[X_1,\dots,X_{n-1},Y_i],\dots,Y_n]~.
\end{equation}
The equation \eqref{eq:n-Jacobi} turns $\VV$ into a \emph{$n$-Lie
  algebra}, a notion introduced in \cite{Filippov} and studied since
by many authors.\footnote{This structure is sometimes also called a
  \emph{Filippov} algebra.} (Notice that, perhaps unfortunately, in
this notation, a Lie algebra is a $2$-Lie algebra.)  More is true,
however, and again the fact that $F\in\Lambda^{n+1}\VV^*$ means that
\begin{equation}
  \label{eq:metrick-Lie}
  \left<[X_1,\dots,X_{n-1},X_n],X_{n+1}\right> = -
  \left<[X_1,\dots,X_{n-1},X_{n+1}],X_n\right>~,
\end{equation}
which we tentatively call a $n$-Lie algebra with an invariant metric,
or a \emph{metric $n$-Lie algebra} for short.

To see that equations \eqref{eq:conj} and \eqref{eq:n-Jacobi} are the
same, let us first rewrite equation \eqref{eq:conj} as follows:
\begin{equation*}
  \sum_a \iota_{\boldsymbol{X}} F^a \wedge F_a = 0~,
\end{equation*}
where $\boldsymbol{X}$ stands for a ($n-1$)-vector
$X_1\wedge\cdots\wedge X_{n-1}$, and where $F_a = \iota_{e_a} F$ and
$F^a=\iota_{e^a} F$ with $e_a = g_{ab} e^b$.  Contracting the above
equation with $n+1$ vectors $Y_1$,...,$Y_{n+1}$, we obtain
\begin{equation*}
  \sum_a \left( \iota_{\boldsymbol{X}} F^a \wedge F_a\right)
  (Y_1,Y_2,\dots,Y_{n+1}) = 0~,
\end{equation*}
which can be rewritten as
\begin{equation*}
  \sum_{i=1}^{n+1} (-1)^{i-1} \left<
  [X_1,\dots,X_{n-1},Y_i], [
  Y_1,\dots,\widehat{Y_i},\dots,Y_{n+1}]\right> = 0~,
\end{equation*}
where a hat over a symbol denotes its omission.  This equation is
equivalent to
\begin{multline*}
  \left<[X_1,\dots,X_{n-1},Y_{n+1}], [Y_1,\dots,Y_n]\right>\\
  = \sum_{i=1}^n (-1)^{n-i} \left< [X_1,\dots,X_{n-1},Y_i],
    [Y_1,\dots,\widehat{Y_i},\dots,Y_{n+1}]\right>~.
\end{multline*}
Finally we use the invariance property \eqref{eq:metrick-Lie} of the
metric to arrive at
\begin{multline*}
  \left<[X_1,\dots,X_{n-1},[Y_1,\dots,Y_n]], Y_{n+1}\right>\\
  = \sum_{i=1}^n \left< [Y_1,\dots,[X_1,\dots,X_{n-1},Y_i],\dots,Y_n],
  Y_{n+1}\right>~,
\end{multline*}
which, since this is true in particular for all $Y_{n+1}$, agrees with
\eqref{eq:n-Lie}.

There seems to be some structure theory for $n$-Lie algebras but to
our knowledge so far nothing on metric $n$-Lie algebras.  Developing
this theory further one could perhaps gain further insight into this
conjecture.  We are not aware of a notion of $n$-Lie group, but if it
did exist then both $\AdS_5 \times S^5$ and the IIB Hpp-wave would be
examples of $4$-Lie groups!

\section{Verifications in low dimension}

To verify the conjecture in the cases mentioned above, we shall use
some group theory and the fact that any two-form can be
skew-diagonalised by an orthogonal transformation, to write down an
ansatz for the $p$-form which we then proceed to analyse
systematically.  Some of the calculations leading to the verification
of the conjecture have been done or checked with \textsl{Mathematica}
and are contained in notebooks which are available upon request.
Since the inner product allows us to identify $\VV$ and its dual
$\VV^*$, we will ignore the distinction in what follows.

\subsection{Proof for $F \in \Lambda^3\EE^6$}

Let $F \in \Lambda^3 \EE^6$ be a $3$-form in six-dimensional euclidean
space.  There is an orthonormal basis $\{e_1,e_2,\ldots,e_6\}$ for
which the $2$-form $\iota_1 F$ obtained by contracting $e_1$ into $F$
takes the form
\begin{equation*}
  \iota_1 F = \alpha e_{23} + \beta e_{45}~,
\end{equation*}
where $e_{ij} = e_i \wedge e_j$ and similarly for $e_{ij\dots k}$ in
what follows.

We must distinguish several cases depending on whether $\alpha$ and
$\beta$ are generic or not.  In the general case, $\iota_1 F$ is a
generic element of a Cartan subalgebra of $\fso(4)$ acting on $\EE^4
=\RR\left<e_2,e_3,e_4,e_5\right>$.  The non-generic cases are in
one-to-one correspondence with conjugacy classes of subalgebras of
$\fso(4)$ of strictly lower rank.  In summary we have the following
cases to consider:
\begin{enumerate}
\item $\fso(4)$: $\alpha$ and $\beta$ generic,
\item $\fsu(2)$: $\alpha = \pm\beta \neq 0$, and
\item $\fso(2)$: $\beta = 0$, $\alpha \neq 0$.
\end{enumerate}
We now treat each case in turn.

\subsubsection{$\fso(4)$}

In the first case, $\alpha$ and $\beta$ are generic, whence the
equation $[\iota_1 F, F] = 0$ says that only terms invariant under the
maximal torus generated by $\iota_1 F$ survive, whence
\begin{equation*}
  F = \alpha e_{123} + \beta e_{145} + \gamma e_{236} + \delta
  e_{456}~.
\end{equation*}
The remaining equations $[\iota_i F, F] = 0$ are satisfied if and only
if
\begin{equation}
  \label{eq:e3so4}
  \alpha \beta + \gamma \delta = 0~.
\end{equation}
Therefore we see that indeed
\begin{equation*}
  F = (\alpha e_1 + \gamma e_6) \wedge e_{23} + (\beta e_1 + \delta
  e_6)\wedge e_{45}
\end{equation*}
can be written as the sum of two simple forms which moreover are
orthogonal, since equation \eqref{eq:e3so4} implies that
\begin{equation*}
  (\alpha e_1 + \gamma e_6)  \perp  (\beta e_1 + \delta e_6)~.
\end{equation*}

\subsubsection{$\fsu(2)$}

Suppose that $\alpha = \beta$ (the case $\alpha = -\beta$ is similar),
so that
\begin{equation*}
  \iota_1 F = \alpha (e_{23} + e_{45})~.
\end{equation*}
This means that $\iota_1 F$ belongs to the Cartan subalgebra of the
selfdual $\SU(2)$ in $\SO(4)$.  The condition $[\iota_1 F, F] = 0$
implies that only terms which have zero weights with respect to this
selfdual $\fsu(2)$ survive, whence
\begin{multline*}
  F = \alpha (e_{123} + e_{145})\\
  + e_6 \wedge \left( \eta(e_{23} + e_{45}) + \gamma(e_{23} - e_{45})
    + \delta (e_{34} - e_{25}) + \varepsilon (e_{24} +
    e_{35})\right)~.
\end{multline*}
However we are allowed to rotate the basis by the normaliser of this
Cartan subalgebra, which is $\U(1) \times \SU(2)$, where the $\U(1)$
is the circle generated by $\iota_1 F$ and the $\SU(2)$ is
anti-selfdual.  Conjugating by the anti-selfdual $\SU(2)$ means that
we can put $\delta = \varepsilon = 0$, say.  The remaining equations
$[\iota_X F,F] =0$ are satisfied if and only if
\begin{equation}
  \label{eq:e3su2}
  \alpha^2 + \eta^2 = \gamma^2~.
\end{equation}
This means that
\begin{equation*}
  F = (\alpha e_1 + (\eta + \gamma)e_6) \wedge e_{23} + (\alpha e_1 +
  (\eta - \gamma) e_6) \wedge e_{45}~,
\end{equation*}
whence $F$ can indeed be written as a sum of two simple $3$-form which
moreover are orthogonal since equation \eqref{eq:e3su2}
implies that
\begin{equation*}
  (\alpha e_1 + (\eta + \gamma)e_6) \perp (\alpha e_1 +
  (\eta - \gamma) e_6)~,
\end{equation*}
as desired.

\subsubsection{$\fso(2)$}

Finally let us consider the case where
\begin{equation*}
  \iota_1 F = \alpha e_{23}~.
\end{equation*}
The surviving terms in $F$ after applying $[\iota_1 F, F] = 0$, are
\begin{equation*}
  F = \alpha e_{123} + \eta e_{234} + \gamma e_{235} + \delta e_{236}
  + \varepsilon e_{456}~.
\end{equation*}
But we can rotate in the (456) plane to make $\gamma = \delta = 0$,
whence
\begin{equation*}
    F = (\alpha e_1 + \eta e_4) \wedge e_{23} + \varepsilon e_4
    \wedge e_{56}
\end{equation*}
can be written as a sum of two simple forms.  Finally the remaining
equations $[\iota_X F,F]=0$ simply say that
\begin{equation}
  \label{eq:e3so2}
  \eta\varepsilon = 0~,
\end{equation}
whence the simple forms are orthogonal, since
\begin{equation*}
  (\alpha e_1 + \eta e_4) \perp \varepsilon e_4~.
\end{equation*}
This verifies the conjecture for $d=3$ and euclidean signature.

\subsection{Proof for $F \in \Lambda^3\EE^{1,5}$}

The lorentzian case is almost identical to the euclidean case, with a
few signs in the equations distinguishing them.  Let
$F\in\Lambda^3\EE^{1,5}$ be a $3$-form in six-dimensional Minkowski
spacetime with pseudo-orthonormal basis $\{e_0,e_2,\ldots,e_6\}$ with
$e_0$ timelike.  Rotating if necessary in the five-dimensional
euclidean space spanned by $\{e_2,e_3,\ldots,e_6\}$, we can guarantee
that
\begin{equation*}
  \iota_0 F = \alpha e_{23} + \beta e_{45}~,
\end{equation*}
as for the euclidean case.  As in that case, we must distinguish
between three cases:
\begin{enumerate}
\item $\fso(4)$: $\alpha$ and $\beta$ generic,
\item $\fsu(2)$: $\alpha = \pm\beta \neq 0$, and
\item $\fso(2)$: $\beta = 0$, $\alpha \neq 0$,
\end{enumerate}
which we now briefly treat in turn.

In the first case, $[\iota_0 F, F]= 0$ means that the only terms in
$F$ which survive are
\begin{equation*}
  F = \alpha e_{023} + \beta e_{045} + \gamma e_{236} + \delta
  e_{456}~,
\end{equation*}
which is already a sum of two simple forms
\begin{equation*}
  F = (\alpha e_0 + \gamma e_6) \wedge e_{23} + ( \beta e_0 + \delta
  e_6) \wedge e_{45}~.
\end{equation*}
The remaining equations $[\iota_X F, F] = 0$ are satisfied if and only
if
\begin{equation}
  \label{eq:m3so4}
  \alpha\beta = \gamma\delta~,
\end{equation}
which makes $\alpha e_0 + \gamma e_6$ and  $\beta e_0 + \delta e_6$
orthogonal, verifying the conjecture in this case. We remark that
this includes the null case as stated in \cite{FOPMax} which
corresponds to setting $\alpha=\beta=\gamma=\delta$.

In the second case, let $\iota_0 F = \alpha (e_{23} + e_{45})$, with
the other possibility $\alpha = -\beta$ being similar.  The equation
$[\iota_0 F, F] = 0$ results in the following
\begin{multline*}
  F = \alpha (e_{023} + e_{045})\\ + e_6 \wedge \left( \eta (e_{23} +
    e_{45}) + \gamma(e_{23} - e_{45}) + \delta (e_{24} + e_{35}) +
    \varepsilon (e_{25} +  e_{34})\right)~.
\end{multline*}
We can rotate by the anti-selfdual $\SU(2) \subset \SO(4)$ in such a
way that $\delta = \varepsilon = 0$, whence $F$ take the desired form
\begin{equation*}
  F = (\alpha e_0 + (\eta + \gamma) e_6) \wedge e_{23} + (\alpha e_0
  + (\eta-\gamma) e_6) \wedge e_{45}~.
\end{equation*}
The remaining equations $[\iota_X F, F] = 0$ are satisfied if and only
if
\begin{equation}
  \label{eq:m3su2}
  \alpha^2 + \gamma^2  = \eta^2~,
\end{equation}
which makes $\alpha e_0 + (\eta + \gamma) e_6$ and $\alpha e_0 +
(\eta - \gamma) e_6$ orthogonal, verifying the conjecture in this
case.

Finally let $\iota_0 F = \alpha e_{23}$.  The equation $[\iota_0 F, F]
= 0$ implies that
\begin{equation*}
  F = \alpha e_{023} + \eta e_{234} + \gamma e_{235} + \delta e_{236}
  + \varepsilon e_{456}~.
\end{equation*}
Rotating in the (456) plane we can make $\gamma = \delta = 0$, whence
$F$ takes the desired form
\begin{equation*}
  F = (\alpha e_0 + \eta e_4) \wedge e_{23} + \varepsilon e_4 \wedge
  e_{56}~.
\end{equation*}
The remaining equations $[\iota_X F, F] = 0$ are satisfied if and only
if
\begin{equation}
  \label{eq:m3so2}
  \eta \varepsilon = 0~,
\end{equation}
making $\alpha e_0 + \eta e_4$ and $\varepsilon e_4$ orthogonal, and
verifying the conjecture in this case, and hence in general for $d=3$
and lorentzian signature.

\subsection{Proof for $F \in \Lambda^3\EE^7$}

Let $F \in \Lambda^3 \EE^7$ be a three-form in a seven-dimensional
euclidean space with orthonormal basis $\{e_i\}_{i=1,\dots,7}$,
relative to which the $2$-form $\iota_7 F$ obtained by contracting
$e_7$ into $F$ takes the form
\begin{equation*}
  \iota_7 F = \alpha e_{12} + \beta e_{34}+ \gamma e_{56}~,
\end{equation*}
where $e_{ij} = e_i \wedge e_j$ and similarly for $e_{ij\dots k}$ in
what follows.

We must distinguish several cases depending on whether $\alpha$,
$\beta$ and $\gamma$ are generic or not.  In the general case,
$\iota_7 F$ is a generic element of  a Cartan subalgebra of $\fso(6)$
acting on the euclidean space $\EE^6$ spanned by
$\{e_i\}_{i=1,\dots,6}$.  The non-generic cases are in one-to-one
correspondence with conjugacy classes of subalgebras of $\fso(6)$ of
strictly lower rank.  In summary we have the following cases to
consider:
\begin{enumerate}
\item $\fso(6)$: $\alpha$, $\beta$ and $\gamma$ generic;
\item $\fsu(2)\times \fu(1)$: $\alpha =\pm \beta$ and $\gamma$ generic;
\item $\fu(1)$ diagonal: $\alpha = \beta= \gamma$;
\item $\fsu(3)$: $\alpha + \beta + \gamma=0$;
\item $\fso(4)$: $\alpha$, $\beta$ generic and $\gamma=0$;
\item $\fsu(2)$: $\alpha=\pm\beta$ and $\gamma=0$; and
\item $\fso(2)$: $\gamma=\beta = 0$, $\alpha \neq 0$.
\end{enumerate}
We now treat each case in turn.

\subsubsection{$\fso(6)$}

In the first case, $\alpha$, $\beta$ and $\gamma$ are generic, whence
the equation $[\iota_7 F, F] = 0$ says that only terms invariant under
the maximal torus generated by $\iota_7 F$ survive, whence
\begin{equation*}
  F = \alpha e_{127} + \beta e_{347} + \gamma e_{567} ~.
\end{equation*}
The remaining equations $[\iota_i F, F] = 0$ are satisfied if and only
if two of $\alpha$, $\beta$ and $\gamma$ vanish, violating the
hypothesis.

\subsubsection{$\fsu(2)\times \fu(1)$}

We choose $\beta=\gamma$ and $\alpha$ generic. The case
$\beta=-\gamma$ is similar.  The equation $[\iota_7 F, F] = 0$ says
that only terms invariant under the maximal torus generated by
$\iota_7 F$ survive. Thus
\begin{multline*}
  F =\alpha e_{127}+ \beta \left(e_{347} + e_{567}\right)\\
  + e_7 \wedge \left(  \delta(e_{34} - e_{56})
    + \varepsilon (e_{36} - e_{45}) + \eta (e_{25} +
    e_{46})\right)~.
\end{multline*}
Using an anti-selfdual rotation, we can set $\varepsilon=\eta=0$. If
$\delta\not=0$, then $\beta+\delta\not=\beta-\delta$ and this leads to
the case investigated in the previous section. If $\delta=0$,
invariance under $[\iota_1 F, F] = 0$ implies that either $\alpha$ or
$\beta$ vanishes, which violates the hypothesis.

\subsubsection{$\fu(1)$ diagonal}
Suppose that $\alpha=\beta=\gamma$. The equation $[\iota_7 F, F] = 0$
implies that
\begin{equation*}
  F =\alpha\left(e_{127}+ e_{347} + e_{567}\right)~.
\end{equation*}
In addition invariance under $[\iota_1 F, F] = 0$ implies that
$\alpha=0$ which violates the hypothesis.

\subsubsection{$\fsu(3)$}
Suppose that $\alpha+\beta+\gamma=0$.  The condition $[\iota_7 F, F] =
0$ implies that
\begin{equation*}
  F =\left(\alpha e_{127}+\beta e_{347} +\gamma e_{567}\right) +
  \delta\Omega_1+ \varepsilon \Omega_2~,
\end{equation*}
where $\Omega_1$ and the real and imaginary parts of the
$\fsu(3)$-invariant (3,0)-form with respect to a complex structure
$J=e_{12}+e_{34}+e_{56}$, that is,
\begin{equation}
  \label{eq:hol3form1}
  \begin{aligned}[m]
  \Omega_1&=e_{135}-e_{146}-e_{236}-e_{245}\\
   \Omega_2&= e_{136}+e_{145}+e_{235}-e_{246}~.
   \end{aligned}
\end{equation}
The presence of these forms can be seen from the decomposition of
$\Lambda^3 \EE^6$ representation under $\fsu(3)$.  Under $\fsu(3)$,
the representation $\EE^6$ transforms as the underlying real
representation of $\repre{3} \oplus \repre{\bar 3}$ (or
$[\![\repre{3}]\!]$ in Salamon's notation \cite{Salamon}).  Similarly
the representation $\Lambda^3 \EE^6$ decomposes into
\begin{equation*}
    \Lambda^3 \EE^6 = [\![ \repre{1}]\!] \oplus [\![\repre{6}]\!] \oplus
    [\![\repre{3}]\!]~.
\end{equation*}
The invariant forms are associated with the trivial representations in
the decomposition.  We still have the freedom to rotate by the
normaliser in $SO(6)$ of the maximal torus of $SU(3)$. An obvious
choice is the diagonal $U(1)$ subgroup of $U(3)$ which leaves
invariant $J$. This $U(1)$ rotates $\Omega_1$ and $\Omega_2$ and we
can use it to set $\varepsilon=0$.  The new case is when
$\delta\not=0$. In such case invariance under the rest of the rotation
$\iota_iF$ implies that $\alpha\beta+2\delta^2=0$ and cyclic in
$\alpha$, $\beta$ and $\gamma$. These relations contradict the
hypothesis that $\alpha+\beta+\gamma=0$ but otherwise generic.

\subsubsection{$\fso(4)$}

Suppose that $\alpha$ and $\beta$ are generic and $\gamma=0$. In that
case, $[\iota_7 F, F] = 0$ implies that
\begin{equation*}
  F =\alpha e_{127}+\beta e_{347}
  +\delta_1 e_{125}+\delta_2 e_{126}
  + \varepsilon_1 e_{345}+\varepsilon_2 e_{346}~.
\end{equation*}
Using a rotation in the (56) plane, we can set $\delta_2=0$. In
addition $\delta_1$ can also be set to zero with a rotation in the
(57) plane and appropriate redefinition of the $\alpha$, $\beta$ and
$\varepsilon_1$ components. Thus the three-form can be written as
\begin{equation*}
  F =\alpha e_{127}+\beta e_{347}
  + \varepsilon_1 e_{345}+\varepsilon_2 e_{346}~.
\end{equation*}
A rotation in the (56) plane leads to $\varepsilon_2=0$. The rest of
the conditions $[\iota_i F, F] = 0$ imply that $\alpha\beta=0$ which
proves the conjecture.

\subsubsection{$\fsu(2)$}
Suppose that $\alpha=\beta$ and $\gamma=0$. The case $\alpha=-\beta$
can be treated similarly. The condition $[\iota_7 F, F] = 0$ implies
that
\begin{multline*}
  F =\alpha \left(e_{127}+e_{347}\right)
  +\delta \left(e_{125}+e_{345}\right)+\varepsilon
  \left(e_{126}+e_{346}\right)\\
  +\eta_1 \left(e_{125}-e_{345}\right)+\eta_2
  \left(e_{145}-e_{235}\right) + \eta_3 \left(e_{135}+e_{245}\right)\\
  +\theta_1\left(e_{126}-e_{346}\right)+\theta_2 \left(e_{146}-e_{236}\right)
  +\theta_3 \left(e_{136}+e_{246}\right)~.
\end{multline*}
With an anti-selfdual rotation, we can set $\eta_2=\eta_3=0$. There
are two cases to consider. If $\eta_1\not=0$, the condition $[\iota_5
F, F] = 0$ implies that $\theta_2=\theta_3=0$. In such case $F$ can be
rewritten as:
\begin{multline*}
  F =\left(\alpha e_7 + (\delta + \eta_1)
    e_5 + (\varepsilon + \theta_1) e_6 \right)\wedge e_{12}\\
  +\left(\alpha e_7 + (\delta - \eta_1) e_5 + (\varepsilon-\theta_1)
    e_6 \right)\wedge e_{34}~.
\end{multline*}
The rest of the conditions imply that
\begin{equation*}
  \alpha^2+\delta^2-\eta_1^2+\epsilon^2-\theta_1^2=0
\end{equation*}
and so $F$ is the sum of two orthogonal simple forms.

Now if $\eta_1=0$, an anti-selfdual rotation will give
$\theta_2=\theta_3=0$.  This case is a special case of the previous
one for which $\eta_1=0$. The conjecture is confirmed.

\subsubsection{$\fso(2)$}

Suppose that $\alpha\not=0$ and $\beta=\gamma=0$.  The condition
$[\iota_7 F, F] = 0$ implies that
\begin{multline*}
  F =\alpha e_{127}+ \sigma_1 e_{123}+\sigma_2 e_{124}+\sigma_3
  e_{125} + \sigma_4 e_{126}\\
  +\tau_1 e_{345}+\tau_2 e_{346}+\tau_3 e_{456}~.
\end{multline*}
A rotation in the (3456) plane can lead to
$\sigma_2=\sigma_3=\sigma_4=0$.  If $\sigma_1\not=0$, then the
condition $[\iota_1 F, F] = 0$ implies that $\tau_2=\tau_1=0$ in which
case
\begin{equation*}
  F =\alpha e_{127}+ \sigma_1 e_{123}+\tau_3 e_{456}~.
\end{equation*}
A further rotation in the $37$ plane leads to the desired result.

Now if $\sigma_1=0$, a rotation in the (3456) plane can lead to
$\tau_2=\tau_3=0$ in which case
\begin{equation*}
  F =\alpha e_{127} + \tau_1 e_{345}~.
\end{equation*}
This again gives the desired result.

\subsection{Proof for $F \in \Lambda^3\EE^d$ and $F \in \Lambda^3\EE^{1,d-1}$, $d<6$}

We shall focus on the proof of the conjecture for $F \in \Lambda^3\EE^d$. The proof of the
statement in the lorentzian case is similar.
Let $F \in \Lambda^3 \EE^5$ be a $3$-form in five-dimensional euclidean
space.  There is an orthonormal basis $\{e_1,e_2,\ldots,e_5\}$ for
which  $\iota_1 F$
takes the form
\begin{equation*}
  \iota_1 F = \alpha e_{23} + \beta e_{45}~.
\end{equation*}

As previous cases, there are  several possibilities to consider depending on whether $\alpha$ and
$\beta$ are generic or not.   Using the adopted group theoretic characterization,
we have the following
cases:
\begin{enumerate}
\item $\fso(4)$: $\alpha$ and $\beta$ generic,
\item $\fsu(2)$: $\alpha = \pm\beta \neq 0$, and
\item $\fso(2)$: $\beta = 0$, $\alpha \neq 0$.
\end{enumerate}
We now treat each case in turn.

\subsubsection{$\fso(4)$}

In the first case, $\alpha$ and $\beta$ are generic, whence the
equation $[\iota_1 F, F] = 0$ says that only terms invariant under the
maximal torus generated by $\iota_1 F$ survive, whence
\begin{equation*}
  F = \alpha e_{123} + \beta e_{145}
\end{equation*}
The remaining equations $[\iota_i F, F] = 0$ are satisfied if and only
if
\begin{equation}
  \label{eq:e3so45}
  \alpha \beta  = 0~,
\end{equation}
which is a contradiction. Thus $\iota_1 F$ cannot be generic.

\subsubsection{$\fsu(2)$}

Suppose that $\alpha = \beta$ (the case $\alpha = -\beta$ is similar),
so that
\begin{equation*}
  \iota_1 F = \alpha (e_{23} + e_{45})~.
\end{equation*}
This means that $\iota_1 F$ belongs to the Cartan subalgebra of the
selfdual $\SU(2)$ in $\SO(4)$.  The condition $[\iota_1 F, F] = 0$
implies that only terms which have zero weights with respect to this
selfdual $\fsu(2)$ survive, and so
\begin{equation*}
  F = \alpha (e_{123} + e_{145})~.
\end{equation*}
 The remaining equations
$[\iota_X F,F] =0$ are satisfied if and only if
\begin{equation}
  \label{eq:e3su25}
  \alpha^2  = 0~,
\end{equation}
which is a contradiction. Thus $\iota_1 F$ cannot be selfdual.

\subsubsection{$\fso(2)$}

Finally let us consider the case where
\begin{equation*}
  \iota_1 F = \alpha e_{23}~.
\end{equation*}
The surviving terms in $F$ after applying $[\iota_1 F, F] = 0$, are
\begin{equation*}
  F = \alpha e_{123} + \eta e_{234} + \gamma e_{235} ~.
\end{equation*}
But we can rotate in the (45) plane to make $\gamma = 0$,
whence
\begin{equation*}
    F = (\alpha e_1 + \eta e_4) \wedge e_{23}
\end{equation*}
is a simple form. This verifies the
conjecture for  $d=5$ and euclidean signature.

\subsubsection{Proof for $F \in \Lambda^3\EE^d$ and $F \in \Lambda^3\EE^{1,d-1}$, $d=3,4$}

The proof for $d=3$ is obvious. It remains to show the
conjecture for $d=4$. In euclidean signature, we have
\begin{equation*}
  \iota_1 F = \alpha e_{23}~.
\end{equation*}
The surviving terms in $F$ after applying $[\iota_1 F, F] = 0$, are
\begin{equation*}
  F = \alpha e_{123} + \eta e_{234} ~.
\end{equation*}
which can be rewritten as
\begin{equation*}
    F = (\alpha e_1 + \eta e_4) \wedge e_{23}
\end{equation*}
and so it is a simple form. This verifies the
conjecture for  $d=4$ and euclidean signature.
The proof for lorentzian spaces is similar.

\subsection{Metric Lie algebras and the case $p{=}3$}

We can give an alternate proof for the case $p{=}3$ exploiting the
relationship with metric Lie algebras; that is, Lie algebras admitting
an invariant non-degenerate scalar product.

It is well-known that reductive Lie algebras --- that is, direct
products of semisimple and abelian Lie algebras --- admit invariant
scalar products: Cartan's criterion allows us to use the Killing form
on the semisimple factor and any scalar product on an abelian Lie
algebra is automatically invariant.

Another well-known example of Lie algebras admitting an invariant
scalar product are the classical doubles.  Let $\fh$ be \emph{any} Lie
algebra and let $\fh^*$ denote the dual space on which $\fh$ acts via
the coadjoint representation.  The definition of the coadjoint
representation is such that the dual pairing $\fh \otimes \fh^* \to
\RR$ is an invariant scalar product on the semidirect product $\fh
\ltimes \fh^*$ with $\fh^*$ an abelian ideal.  The Lie algebra $\fh
\ltimes \fh^*$ is called the classical double of $\fh$ and the
invariant metric has split signature $(r,r)$ where $\dim \fh = r$.

It turns out that all Lie algebras admitting an invariant scalar
product can be obtained by a mixture of these constructions.  Let
$\fg$ be a Lie algebra with an invariant scalar product
$\left<-,-\right>_\fg$, and let $\fh$ act on $\fg$ preserving both the
Lie bracket and the scalar product; in other words, $\fh$ acts on
$\fg$ via skew-symmetric derivations.  First of all, since $\fh$ acts
on $\fg$ preserving the scalar product, we have a linear map
\begin{equation*}
  \fh \to \fso(\fg) \cong \Lambda^2 \fg~,
\end{equation*}
with dual map
\begin{equation*}
  c: \Lambda^2 \fg^* \cong \Lambda^2 \fg \to \fh^*~,
\end{equation*}
where we have used the invariant scalar product to identity $\fg$ and
$\fg^*$ equivariantly.  Since $\fh$ preserves the Lie bracket in
$\fg$, this map is a cocycle, whence it defines a class $[c]\in
H^2(\fg;\fh^*)$ in the second Lie algebra cohomology of $\fg$ with
coefficients in the trivial module $\fh^*$.  Let $\fg \times_c \fh^*$
denote the corresponding central extension.  The Lie bracket of the
$\fg \times_c \fh^*$ is such that $\fh^*$ is central and if
$X,Y\in\fg$, then
\begin{equation*}
  [X,Y] = [X,Y]_\fg + c(X,Y)~,
\end{equation*}
where $[-,-]_\fg$ is the Lie bracket of $\fg$.  Now $\fh$ acts
naturally on this central extension: the action on $\fh^*$ given by
the coadjoint representation.  This then allows us to define the
\emph{double extension} of $\fg$ by $\fh$,
\begin{equation*}
  \fd(\fg,\fh) = \fh \ltimes (\fg \times_c \fh^*)
\end{equation*}
as a semidirect product.  Details of this construction can be found in
\cite{MedinaRevoy,FSsug}.  The remarkable fact is that $\fd(\fg,\fh)$
admits an invariant inner product:
\begin{equation}
  \bordermatrix{& \fg & \fh & \fh^* \cr
  \fg &  \left<-,-\right>_\fg & 0         &  0        \cr
  \fh &    0   & B & \id \cr
  \fh^*&    0   & \id &  0 \cr}
\end{equation}
where $B$ is \emph{any} invariant symmetric bilinear form on $\fh$ and
$\id$ stands for the dual pairing between $\fh$ and $\fh^*$.

We say that a Lie algebra with an invariant scalar product is
indecomposable if it cannot be written as the direct product of two
orthogonal ideals.  A theorem of Medina and Revoy \cite{MedinaRevoy}
(see also \cite{FSalgebra} for a refinement) says that an
indecomposable (finite-dimensional) Lie algebra with an invariant
scalar product is one of the following:
\begin{enumerate}
\item one-dimensional,
\item simple, or
\item a double extension $\fd(\fg,\fh)$ where $\fh$ is either simple
  or one-dimensional and $\fg$ is a Lie algebra with an invariant
  scalar product.  (Notice that we can take $\fg$ to be the trivial
  zero-dimensional Lie algebra.  In this way we recover the classical
  double.)
\end{enumerate}
Any (finite-dimensional) Lie algebra with an invariant scalar product
is then a direct sum of indecomposables.

Notice that if the scalar product on $\fg$ has signature $(p,q)$ and
if $\dim \fh = r$, then the scalar product on $\fd(\fg,\fh)$ has
signature $(p+r,q+r)$.  Therefore euclidean Lie algebras are
necessarily reductive, and if indecomposable they are either
one-dimensional or simple.  Up to dimension $7$ we have the following
euclidean Lie algebras:
\begin{itemize}
\item $\RR^d$ with $d\leq 7$,
\item $\fsu(2) \oplus \RR^k$ with $k\leq 4$, and
\item $\fsu(2) \oplus \fsu(2) \oplus \RR^k$ with $k=0,1$.
\end{itemize}
The conjecture clearly holds for all of them.

The lorentzian case is more involved.  Indecomposable lorentzian Lie
algebras are either reductive or double extensions $\fd(\fg,\fh)$
where $\fg$ has a positive-definite invariant scalar product and $\fh$
is one-dimensional.  In the reductive case, indecomposability means
that it has to be simple, whereas in the latter case, since the scalar
product on $\fg$ is positive-definite, $\fg$ must be reductive.  A
result of \cite{FSsug} (see also \cite{FSalgebra}) then says that any
semisimple factor in $\fg$ splits off resulting in a decomposable Lie
algebra.  Thus if the double extension is to be indecomposable, then
$\fg$ must be abelian.  In summary, an indecomposable lorentzian Lie
algebra is either simple or a double extension of an abelian Lie
algebra by a one-dimensional Lie algebra and hence solvable (see,
e.g., \cite{MedinaRevoy}).

These considerations make possible the following enumeration of
lorentzian Lie algebras up to dimension $7$:
\begin{enumerate}
\item $\EE^{1,d-1}$ with $d\leq 7$,
\item $\EE^{1,k} \oplus \fso(3)$ with $k\leq 3$,
\item $\EE^k \oplus \fso(1,2)$ with $k \leq 4$,
\item $\fso(1,2) \oplus \fso(3) \oplus \EE^k$ with $k=0,1$, or
\item $\fd(\EE^4,\RR)\oplus \EE^k$ with $k=0,1$,
\end{enumerate}
where the last case actually corresponds to a family of Lie algebras,
depending on the action of $\RR$ on $\EE^4$.  The conjecture holds
manifestly for all cases except possibly the last, which we must
investigate in more detail.

Let $e_i$, $i=1,2,3,4$, be an orthonormal basis for $\EE^4$, and let
$e_- \in \RR$ and $e_+\in\RR^*$, so that together they span
$\fd(\EE^4,\RR)$.  The action of $\RR$ on $\RR^4$ defines a map $\rho:
\RR \to \Lambda^2 \RR^4$, which can be brought to the form $\rho(e_-) =
\alpha e_1 \wedge e_2 + \beta e_3 \wedge e_4$ via an orthogonal change
of basis in $\EE^4$ which moreover preserves the orientation.  The Lie
brackets of $\fd(\EE^4,\RR)$ are given by
\begin{equation*}
  \begin{aligned}[m]
    [e_-,e_1] &= \alpha e_2\\
    [e_-,e_2] &= -\alpha e_1\\
    [e_1,e_2] &= \alpha e_+
  \end{aligned}\qquad\qquad
  \begin{aligned}[m]
    [e_-,e_3] &= \beta e_4\\
    [e_-,e_4] &= -\beta e_3\\
    [e_3,e_4] &= \beta e_+
  \end{aligned}~,
\end{equation*}
and the scalar product is given (up to scale) by
\begin{equation*}
  \left<e_-,e_-\right> = b \qquad \left<e_+,e_-\right> = 1 \qquad
  \left<e_i,e_j\right> = \delta_{ij}~.
\end{equation*}
The first thing we notice is that we can set $b=0$ without loss of
generality by the automorphism fixing all $e_i,e_+$ and mapping $e_-
\mapsto e_- - \half b e_+$.  We will assume that this has been done
and that $\left<e_-,e_-\right>=0$.  A straightforward calculation
shows that the three-form $F$ takes the form
\begin{equation*}
  F = \alpha e_- \wedge e_1 \wedge e_2 + \beta e_- \wedge e_3 \wedge
  e_4~,
\end{equation*}
whence the conjecture holds.

\subsection{Proof for $F \in \Lambda^4\EE^8$}

In the absence (to our knowledge) of a structure theorem for metric
$n$-Lie algebras, we will present the verification of the conjecture
in the remaining cases using the ``brute-force'' approach explained
earlier.

Choose an orthonormal basis $\{e_1,e_2,\ldots,e_8\}$ for which
$\iota_{12} F = \alpha e_{34} + \beta e_{56} + \gamma e_{78}$, where
$\iota_{12}$ means the contraction of $F$ by $e_{12}$.

Suppose that $\alpha$, $\beta$ and $\gamma$ are generic.  In this
case, the equation $[\iota_{12} F,F]=0$ says that the only terms in
$F$ which survive are those which are invariant under the maximal
torus of $\SO(6)$, the group of rotations in the six-dimensional
space spanned by $\{e_3,e_4,\ldots,e_8\}$; that is,
\begin{equation*}
  F = \alpha e_{1234} + \beta e_{1256} + \gamma e_{1278} + \delta
  e_{3456} + \varepsilon e_{3478} + \eta e_{5678}~.
\end{equation*}
Now, $\iota_{13} F = -\alpha e_{24}$, whence the equation $[\iota_{13}
F, F] = 0$ implies that $\beta = \gamma = \delta = \varepsilon = 0$,
violating the condition that $\iota_{12} F$ be generic.

In fact, this argument clearly works for $d\geq 4$ so that for $d\geq
4$ we have to deal with non-generic rotations.  Non-generic rotations
correspond to (conjugacy classes of) subalgebras of $\fso(6)$ with
rank strictly less than that of $\fso(6)$:
\begin{enumerate}
\item $\fsu(3)$: $\alpha + \beta + \gamma = 0$ but all $\alpha$,
  $\beta$, and $\gamma$ nonzero;
\item $\fsu(2) \times \fu(1)$: $\alpha = \beta \neq \gamma$, but again
  all nonzero;
\item $\fu(1)$ diagonal: $\alpha= \beta= \gamma \neq 0$;
\item $\fso(4)$: $\gamma=0$ and $\alpha \neq \beta$ nonzero;
\item $\fsu(2)$: $\gamma=0$ and $\alpha = \beta \neq 0$; and
\item $\fso(2)$: $\beta = \gamma = 0$ and $\alpha \neq 0$.
\end{enumerate}
We now go down this list case by case.

\subsubsection{$\fsu(3)$}

When $\iota_{12}F$ is a generic element of the Cartan subalgebra of an
$\fsu(3)$ subalgebra of $\fso(6)$ the only terms in $F$ which satisfy
the equation $[\iota_{12}F , F]= 0$ are those which have zero weights
relative to this Cartan subalgebra.  Let $\EE^6 =\left< e_1 , e_2
\right>^\perp$.  Then $F$ can be written as
\begin{equation*}
  F = e_{12} \wedge \iota_{12} F + G
\end{equation*}
where $G$ is in the kernel of $\iota_{12}$, namely
\begin{equation*}
  G = e_1 \wedge G_1 + e_2 \wedge G_2 + G_3~,
\end{equation*}
where $G_1,G_2 \in \Lambda^3 \EE^6$ and $G_3 \in \Lambda^4 \EE^6$.
We have investigated the  decomposition of $\Lambda^3 \EE^6$ under $\fsu(3)$
in the previous section. The representation $\Lambda^4 \EE^6$
decomposes into
\begin{equation*}
   \Lambda^4 \EE^6 =  \repre{1} \oplus \repre{8} \oplus
    [\![\repre{3}]\!]~,
 \end{equation*}
whence it is clear where the zero weights are: they are one  in
the trivial representation $\repre{1}$ and two in the adjoint
$\repre{8}$.  This means that in this case together with the
zero weights of the $\Lambda^3 \EE^6$ representations a total of seven
terms in $G$:
\begin{equation*}
  \begin{aligned}[m]
    G_1 &= \lambda_1 \Omega_1 + \lambda_2 \Omega_2\\
    G_2 &= \lambda_3 \Omega_1 + \lambda_4 \Omega_2\\
    G_3 &= \mu_1 e_{3456} + \mu_2 e_{3478} + \mu_3 e_{5678}~,
  \end{aligned}
\end{equation*}
where
\begin{equation}
  \label{eq:hol3form}
  \begin{aligned}[m]
    \Omega_1 &= e_{357} - e_{368} - e_{458} - e_{467}\\
    \Omega_2 &= e_{358} + e_{367} + e_{457} - e_{468}
  \end{aligned}
\end{equation}
are the real and imaginary parts, respectively, of the holomorphic
$3$-form in $\EE^6$ thought of as $\CC^3$ with the $\fsu(3)$-invariant
complex structure $J=e_{34} + e_{56} + e_{78}$.  We still have to
freedom to rotate by the normaliser in $\SO(6)$ of the maximal torus
in $\SU(3)$ that $\iota_{12}F$ determines.  An obvious choice is the
$\U(1)$ generated by the complex structure.  This is not in $\SU(3)$
but in $\U(3)$ and has the virtue of acting on $\Omega = \Omega_1 + i
\Omega_2$ by multiplication by a complex phase.  This means that we
can always choose $\Omega$ to be real, thus setting $\lambda_4 = 0$,
say.  Analysing the remaining equations $[\iota_{ij} F, F]=0$ we see
that $\alpha$ and $\beta$ are constrained to $\alpha = \pm \beta$,
violating the hypothesis that they are generic.

\subsubsection{$\fsu(2) \times \fu(1)$}

Let us consider $\alpha = \beta$, the other case being similar, in
fact related by conjugation in $\O(4)$, which is an outer
automorphism.  The equation $[\iota_{12} F, F]=0$ says that the only
terms in $F$ which survive are those corresponding to zero weights of
the $\fsu(2) \times \fu(1)$ subalgebra of $\fso(6)$.  It is easy to
see that $\Lambda^3\EE^6$ has not zero weights, whereas the zero
weights in $\Lambda^4\EE^6$ are the Hodge duals of the following
$2$-forms:
\begin{equation*}
  e_{34}\qquad e_{56} \qquad e_{78} \qquad e_{35} + e_{46} \qquad
  e_{36} - e_{45}~.
\end{equation*}
Conjugating by the anti-selfdual $\SU(2)$ we can set to zero the
coefficients of the last two forms, leaving
\begin{equation*}
  F =   \alpha (e_{1234} + e_{1256}) + \gamma e_{1278} + \mu_1
  e_{3456} + \mu_2 e_{3478} + \mu_3 e_{5678}
\end{equation*}
as the most general solution of $[\iota_{12}F,F]=0$.  Now the equation
$[\iota_{13}F,F]=0$, for example, implies that $\alpha$ must vanish,
violating the hypothesis.  This case is therefore discarded.

\subsubsection{$\fu(1)$ diagonal}

In this case, $\iota_{12} F = \alpha (e_{34} + e_{56} + e_{78})$
belongs to the diagonal $\fu(1)$ which is the centre of $\fu(3)
\subset \fso(6)$, where $\fso(6)$ acts on the $\EE^6$ spanned by
$\{e_i\}_{3\leq i\leq 8}$.  There are no zero weights in $\Lambda^3
\EE^6$, but there are nine in $\Lambda^4 \EE^6$: the Hodge duals of
$\fu(3) \subset \fso(6) \cong \Lambda^2 \EE^6$.  However we are
allowed to conjugate by the normaliser of $\fu(1)$ in $\fso(6)$ which
is $\fu(3)$.  This allows us to conjugate the invariant $2$-forms to
lie in the Cartan subalgebra of $\fu(3)$.  In summary, the solution to
$[\iota_{12}F,F]=0$ can be written in the form
\begin{equation*}
  F = \alpha (e_{1234} + e_{1256} + e_{1278}) + \mu_1 e_{3456} + \mu_2
  e_{3478} + \mu_3 e_{5678}~.
\end{equation*}
Now we consider for example the equation $[\iota_{13}F,F]=0$ and we
see that $\alpha$ must vanish, violating the hypothesis.  Thus this
case is also discarded.

Notice that all the cases where the $2$-form $\iota_{12}F$ has maximal
rank have been discarded, often after a detailed analysis of the
equations.  This should have a simpler explanation.

\subsubsection{$\fso(4)$}

In this case $\iota_{12}F = \alpha e_{34} + \beta e_{56}$ where
$\alpha$ and $\beta$ are generic.  This means that the most general
solution of $[\iota_{12} F, F]=0$ is given by
\begin{equation*}
  F = \alpha e_{1234} + \beta e_{1256} + G~,
\end{equation*}
where $G$ is of the form $e_1 \wedge G_1 + e_2 \wedge G_2 + G_3$,
where $G_1,G_2 \in \Lambda^3\EE^6$ and $G_3 = \Lambda^4 \EE^6$, where
$\EE^6$ is spanned by $\{e_i\}_{3\leq i\leq 8}$, and where the $G_i$
have zero weight with respect to this $\fso(4)$ algebra.  A little
group theory shows that $G_1$ and $G_2$ are linear combinations of the
four monomials $e_{347}, e_{348}, e_{567}, e_{568}$; whereas $G_3$ is
a linear combination of the three monomials $e_{3456}, e_{3478},
e_{5678}$.  We still have the freedom to conjugate by the normaliser
in $\SO(6)$ of the maximal torus generated by $\iota_{12} F$, which
includes the $\SO(2)$ of rotations in the (78) plane.  Doing this we
can set any one of the monomials in $e_1 \wedge G_1$, say $e_{1347}$,
to zero.  In summary, the most general solution of $[\iota_{12}F,
F]=0$ can be put in the following form
\begin{multline*}
  F = \alpha e_{1234} + \beta e_{1256} + \mu_1 e_{3456} +
  \mu_2 e_{3478} + \mu_3 e_{5678} + \lambda_1 e_{1348}\\
  + \lambda_2 e_{1567} + \lambda_3 e_{1568} + \lambda_4 e_{2347} +
  \lambda_5 e_{2348} + \lambda_6 e_{2567} + \lambda_7 e_{2568}~.
\end{multline*}
Analysing the remaining equations $[\iota_{ij} F, F]=0$ we notice that
genericity of $\alpha$ and $\beta$ are violated unless $\mu_1 = 0$ and
$\mu_3 \mu_2 = \alpha \beta$.  Given this we find that the most
general solution is
\begin{multline*}
  F = \alpha e_{1234} + \beta e_{1256} + \mu_3 e_{5678} +
  \mu_2 e_{3478}+ \nu_1 ( \alpha e_{1348} + \mu_3 e_{2567})\\
  + \nu_2 ( \beta e_{1567} - \mu_2 e_{2348}) + \nu_3 ( \beta
  e_{1568} + \mu_2 e_{2347})
\end{multline*}
subject to
\begin{equation}
  \label{eq:e4so4}
  \nu_1 \nu_3 = - 1 \qquad \text{and}\qquad \mu_3 \mu_2 =
  \alpha\beta~.
\end{equation}
These identities are precisely the ones that allow us to rewrite $F$
as a sum of two simple forms
\begin{equation*}
  \begin{aligned}[m]
    F_1 &= (\alpha e_1 - \mu_2 (\nu_3 e_7 - \nu_2 e_8)) \wedge
    (e_2 + \nu_1 e_8) \wedge e_3 \wedge e_4\\
    F_2 &= (\beta e_1 - \mu_3 \nu_1 e_7) \wedge
    (e_2 + \nu_2 e_7 + \nu_3 e_8) \wedge e_5 \wedge e_6~,
  \end{aligned}
\end{equation*}
which moreover are orthogonal.

\subsubsection{$\fsu(2)$}

In this case $\iota_{12}F = \alpha (e_{34} + e_{56})$, where without
loss of generality we can set $\alpha =1$.  This means that the most
general solution of $[\iota_{12} F, F]=0$ is given by
\begin{equation*}
  F = e_{1234} + e_{1256} + e_1 \wedge G_1 + e_2 \wedge G_2 + G_3~,
\end{equation*}
where $G_1,G_2 \in \Lambda^3\EE^6$ and $G_3 = \Lambda^4 \EE^6$, where
$\EE^6$ is spanned by $\{e_i\}_{3\leq i\leq 8}$, and where the $G_i$
have zero weight with respect to this $\fsu(2)$ algebra.  A little
group theory shows that $G_1$ and $G_2$ are linear combinations of the
following eight $3$-forms
\begin{equation*}
  e_{34i} + e_{56i}\qquad e_{34i} - e_{56i}\qquad e_{35i} +
  e_{46i}\qquad e_{36i} - e_{45i}
\end{equation*}
where $i$ can be either 7 or 8; whereas $G_3$ is the Hodge dual
(in $\EE^6$) of a linear combination of
\begin{equation*}
  e_{34} + e_{56}\qquad e_{34} - e_{56}\qquad e_{35} +
  e_{46}\qquad e_{36} - e_{45}~.
\end{equation*}
Using the freedom to conjugate by the normaliser of $\fsu(2)$ in
$\fso(6)$ we can choose basis such that $G_3$ takes the form
\begin{equation*}
  G_3 = \mu_1 e_{3456} + \mu_2 e_{3478} + \mu_3 e_{5678}~.
\end{equation*}
This means that $F$ takes the following form:
\begin{multline*}
  F = e_{1234} + e_{1256} + \mu_1 e_{3456} +
  \mu_2 e_{3478} + \mu_3 e_{5678}\\
  + \lambda_1 e_{1347} + \lambda_2 e_{1348} + \lambda_3 e_{1567} +
  \lambda_4 e_{1568} + \lambda_5 e_{2347} + \lambda_6 e_{2348}
  + \lambda_7 e_{2567} + \lambda_8 e_{2568}\\
  + \sigma_1 (e_{1357} + e_{1467}) + \sigma_2 (e_{1367} - e_{1457})
  + \sigma_3 (e_{1358} + e_{1468}) + \sigma_4 (e_{1368} - e_{1458})\\
  + \sigma_5 (e_{2357} + e_{2467}) + \sigma_6 (e_{2367} - e_{2457})
  + \sigma_7 (e_{2358} + e_{2468}) + \sigma_8 (e_{2368} - e_{2458})~.
\end{multline*}
This still leaves the possibility of rotating, for example, in the
(78) plane and an anti-selfdual rotation in the (3456) plane.
Rotating in the (78) plane allows us to set $\lambda_8=0$, whereas an
anti-selfdual rotation allows us to set $\sigma_8=0$.  Imposing, for
example, the equation $[\iota_{25}F, F]=0$ tells us that $\lambda_1 =
0$, whereas the rest of the equations also say that $\sigma_2=0$.
It follows after a little work that if $\mu_1 \neq 0$ we arrive at a
contradiction, so that we take $\mu_1 = 0$.

We now have to distinguish between two cases, depending on whether or
not $\mu_2$ equals $\mu_3$.  If $\mu_2\neq\mu_3$, then all $\sigma_i
=0$, and moreover $F$ takes the form
\begin{multline*}
  F = e_{1234} + e_{1256} + \mu_2 e_{3478} + \mu_3 e_{5678}
  + \lambda_2 (e_{1348} + \mu_3 e_{2567})\\
  + \lambda_3 (e_{1567} - \mu_2 e_{2348})
  + \lambda_4 (e_{1568} + \mu_2 e_{2347})~,
\end{multline*}
subject to the equations
\begin{equation}
  \label{eq:e4su2-1}
  \lambda_2 \lambda_4 = -1 \qquad \text{and} \qquad \mu_2 \mu_3 = 1~.
\end{equation}
These equations are precisely what is needed to write $F$ as a sum of
two orthogonal simple forms $F= F_1 + F_2$, where
\begin{equation*}
  \begin{aligned}[m]
    F_1 &= (e_1 - \mu_2 (\lambda_4 e_7 -  \lambda_3 e_8)) \wedge (e_2
    + \lambda_2 e_8) \wedge e_3 \wedge e_4\\
    F_2 &= (e_1 - \mu_3 \lambda_2 e_7) \wedge (e_2 + \lambda_3 e_7 +
    \lambda_4 e_8) \wedge e_5 \wedge e_6~.
  \end{aligned}
\end{equation*}

Finally, we consider the case $\mu_2 = \mu_3$, which is inconsistent
unless $\mu_2^2 = 1$.  Then the most general solution takes the form
\begin{multline*}
  F = e_{1234} + e_{1256} + \mu_2(e_{3478} + e_{5678})\\
  + \lambda_2 (e_{1348} + \mu_2 e_{2567}) + \lambda_3 (e_{1567} -
  \mu_2 e_{2348}) + \lambda_4 (e_{1568} + \mu_2 e_{2347})\\
  + \sigma_1 (e_{1357} + e_{1467} + \mu_2 e_{2358} + \mu_2 e_{2468})\\
  + \sigma_3 (e_{1358} + e_{1468} - \mu_2 e_{2357} - \mu_2 e_{2467})\\
  + \sigma_4 (e_{1368} - e_{1458} - \mu_2 e_{2367} + \mu_2 e_{2457})~,
\end{multline*}
subject to the following equations
\begin{equation}
  \label{eq:e4su2'}
  \begin{gathered}
    \lambda_3 \sigma_4 = 0 = \sigma_1 \sigma_4\\
    (\lambda_2 - \lambda_4) \sigma_1 + \lambda_3 \sigma_3 = 0\\
    \sigma_1^2 + \sigma_3^2 + \sigma_4^2 = 1 + \lambda_2 \lambda_4~.
  \end{gathered}
\end{equation}
Let us rewrite $F$ in terms of (anti)selfdual $2$-forms in the (1278)
and (3456) planes:
\begin{footnotesize}
\begin{multline*}
  F = \left[ (e_{12} + \mu_2 e_{78}) + \half \lambda_3 (e_{17} - \mu_2
    e_{28}) + \half(\lambda_2 + \lambda_4)(e_{18} + \mu_2 e_{27})
  \right] \wedge (e_{34} + e_{56})\\
  + (e_{17} + \mu_2 e_{28}) \wedge \left[ \sigma_1 (e_{35} + e_{46}) -
    \half \lambda_3 (e_{34} - e_{56}) \right]\\
  + (e_{18} - \mu_2 e_{27}) \wedge \left[ \sigma_3 (e_{35} + e_{46}) +
    \sigma_4 (e_{36} - e_{45}) + \half (\lambda_2 - \lambda_4) (e_{34}
    - 3_{56})\right]~.
\end{multline*}
\end{footnotesize}
Notice that the first two equations in \eqref{eq:e4su2'} simply say
that the two anti-selfdual $2$-forms
\begin{gather*}
  \sigma_1 (e_{35} + e_{46}) - \half \lambda_3 (e_{34} - e_{56})\\
 \sigma_3 (e_{35} + e_{46}) + \sigma_4 (e_{36} - e_{45}) + \half
 (\lambda_2 - \lambda_4) (e_{34} - 3_{56})
\end{gather*}
are collinear.  Therefore performing an anti-selfdual rotation in the
$(36)-(45)$ direction, we can eliminate the $e_{35} + e_{46}$ and
$e_{36}-e_{45}$ components, effectively setting
$\sigma_1=\sigma_3=\sigma_4=0$.  This reduces the problem to the
previous case, except that now $\mu_2 = \mu_3$.

\subsubsection{$\fso(2)$}

Finally, we consider the case where $\iota_{12} F = \alpha e_{34}$.
The most general $F$ has the form
\begin{equation*}
  F = \alpha e_{1234} + e_1 \wedge G_1 + e_2 \wedge G_2 + G_3~,
\end{equation*}
where $G_1,G_2 \in \Lambda^3\EE^6$ and $G_3 \in \Lambda^4\EE^6$, where
$\EE^6$ is spanned by $\{e_i\}_{3\leq i\leq 8}$.  Such an $F$ will
obey $[\iota_{12} F, F]=0$ if and only if the $G_i$ have zero weights
under the $\fso(2)$ generated by $\iota_{12}F$.  This means that each
of $G_1,G_2$ is a linear combination of the 8 monomials
\begin{equation*}
  e_{345}\quad e_{346}\quad e_{347}\quad e_{348}\quad
  e_{567}\quad e_{568}\quad e_{578}\quad e_{678}
\end{equation*}
Using the freedom to conjugate by the $\SO(4)$ which acts in the
(5678) plane, we can write the most general $G_3$ as a linear
combination of the monomials $e_{5678}, e_{3478}, e_{3456}$.  This
still leaves the possibility of rotating in the (56)- and (78) planes
separately.  Doing so we can set to zero the coefficients of say,
$e_{2568}$ and $e_{2678}$, leaving a total of 17 free parameters
\begin{multline*}
  F = e_{1234} + \mu_1 e_{3456} + \mu_2 e_{3478} + \mu_3 e_{5678}
  + \lambda_1 e_{1347} + \lambda_2 e_{1348}\\
  + \lambda_3 e_{1567} + \lambda_4 e_{1568} + \lambda_5 e_{2347} +
  \lambda_6 e_{2348} + \lambda_7 e_{2567} + \sigma_1 e_{1345}\\
  + \sigma_2 e_{1346} + \sigma_3 e_{1578} + \sigma_4 e_{1678}
  + \sigma_5 e_{2345} + \sigma_6 e_{2346} + \sigma_7 e_{2578}~,
\end{multline*}
and where we have set $\alpha = 1$ without loss of generality.
We now impose the rest of the equations $[\iota_{ij}F ,F]=0$.  We
first observe that if $\mu_1 \neq 0$, then $\mu_2 = \mu_3 = \lambda_i=
\sigma_3 = \sigma_4 = \sigma_7 = 0$, leaving
\begin{equation*}
  F = e_{1234} + \mu_1 e_{3456} + \sigma_1 e_{1345} + \sigma_2
  e_{1346} +  \sigma_5 e_{2345} + \sigma_6 e_{2346}~,
\end{equation*}
subject to
\begin{equation}
  \label{eq:e4u1-1}
  \sigma_1 \sigma_6 - \sigma_2 \sigma_5 = \mu_1~,
\end{equation}
which guarantees that $F$ is actually a simple form
\begin{equation*}
  F = (e_1 - \sigma_5 e_5 - \sigma_6 e_6) \wedge (e_2 + \sigma_1 e_5 +
  \sigma_2 e_6) \wedge e_3 \wedge e_4~,
\end{equation*}
which is a degenerate case of the conclusion of the conjecture.

Let us then suppose that $\mu_1 = 0$.  We next observe that if $\mu_2
\neq 0$ then $\mu_3 = \sigma_i = \lambda_3 = \lambda_4 = \lambda_7 =
0$.  This is again, up to a relabelling of the coordinates, the same
degenerate case as before and the conclusion still holds.

Finally let us suppose that both $\mu_1$ and $\mu_2$ vanish.  We must
distinguish between two cases, depending on whether $\mu_3$ also
vanishes or not.  If $\mu_3=0$ then we have that $F$ is given by
\begin{multline*}
  F = e_{1234} + \lambda_1 e_{1347} + \lambda_2 e_{1348}
  + \lambda_5 e_{2347} + \lambda_6 e_{2348}\\
  + \sigma_1 e_{1345} + \sigma_2 e_{1346}  + \sigma_5 e_{2345}  +
    \sigma_6 e_{2346}~,
\end{multline*}
subject to the equations
\begin{equation}
  \label{eq:e4u1-2}
  \begin{aligned}[m]
    \lambda_2 \lambda_5 &= \lambda_1 \lambda_6\\
    \lambda_1 \sigma_5 &= \lambda_5 \sigma_1\\
    \lambda_1 \sigma_6 &= \lambda_5 \sigma_2
  \end{aligned}\qquad\qquad\qquad
  \begin{aligned}[m]
    \sigma_2 \sigma_5 &= \sigma_1 \sigma_6\\
    \lambda_6 \sigma_2 &= \lambda_2 \sigma_6\\
    \lambda_6 \sigma_1 &= \lambda_2 \sigma_5~,
  \end{aligned}
\end{equation}
which are precisely the equations which allow us to rewrite $F$ as a
simple form $F = \theta_1 \wedge \theta_2 \wedge e_3 \wedge e_4$,
where
\begin{equation*}
  \begin{aligned}[m]
    \theta_1 &= e_1 - \sigma_5 e_5 - \sigma_6 e_6 - \lambda_5 e_7 -
    \lambda_6 e_8\\
    \theta_2 &= e_2 + \sigma_1 e_5 + \sigma_2 e_6 + \lambda_1 e_7 +
    \lambda_2 e_8~.
  \end{aligned}
\end{equation*}

Finally suppose that $\mu_3 \neq 0$.  In this case $F$ is given by
\begin{multline*}
  F = e_{1234} + \mu_3 e_{5678}
  + \lambda_2 (e_{1348} + \mu_3 e_{2567})\\
  + \lambda_5 (e_{2347} + \mu_3 e_{1568})
  + \lambda_6 (e_{2348} - \mu_3 e_{1567})
  + \sigma_2 (e_{1346} + \mu_3 e_{2578})\\
  + \sigma_5 (e_{2345} + \mu_3 e_{1678})
  + \sigma_6 (e_{2346} - \mu_3 e_{1578})~,
\end{multline*}
subject to the equations
\begin{equation}
  \label{eq:e4u1-3}
  \lambda_2 \lambda_5 = \lambda_2 \sigma_5 = \sigma_2 \lambda_5 =
  \sigma_2 \sigma_5 = 0 \qquad\text{and}\qquad \lambda_6 \sigma_2 =
  \lambda_2 \sigma_6~.
\end{equation}
We must distinguish between three cases:
\begin{enumerate}
\item $\lambda_2 \neq 0$,
\item $\lambda_2 = 0$ and $\sigma_2 \neq 0$, and
\item $\lambda_2 = \sigma_2 = 0$.
\end{enumerate}
We now do each in turn.

If $\lambda_2 \neq 0$, $F$ is given by
\begin{multline*}
  F = e_{1234} + \mu_3 e_{5678}
  + \lambda_2 (e_{1348} + \mu_3 e_{2567})\\
  + \lambda_6 (e_{2348} - \mu_3 e_{1567})
  + \sigma_2 (e_{1346} + \mu_3 e_{2578})
  + \sigma_6 (e_{2346} - \mu_3 e_{1578})~,
\end{multline*}
subject to the second equation in \eqref{eq:e4u1-3}.  This is
precisely the equation that allows us to write $F$ as a sum of two
simple forms $F= F_1 + \mu_3 F_2$, where
\begin{equation*}
  \begin{aligned}[m]
    F_1 &= (e_1 - \sigma_6 e_6 - \lambda_6 e_8) \wedge (e_2 +
    \sigma_2 e_6 + \lambda_2 e_8) \wedge e_3 \wedge e_4\\
    F_2 &= e_5 \wedge (e_6 + \sigma_6 e_1 - \sigma_2 e_2) \wedge e_7
    \wedge (e_8 + \lambda_6 e_1 - \lambda_2 e_2)~.
  \end{aligned}
\end{equation*}
Notice moreover that $F_1$ and $F_2$ are orthogonal.

If $\lambda_2 = 0$ and $\sigma_2 \neq 0$, $F$ is given by
\begin{equation*}
  F = e_{1234} + \mu_3 e_{5678}
  + \sigma_2 (e_{1346} + \mu_3 e_{2578})
  + \sigma_6 (e_{2346} - \mu_3 e_{1578})~,
\end{equation*}
which can be written as a sum $F = F_1 + \mu_3 F_2$ of two simple
forms
\begin{equation*}
  \begin{aligned}[m]
    F_1 &= (e_1 - \sigma_6 e_6) \wedge (e_2 + \sigma_2 e_6) \wedge
    e_3 \wedge e_4\\
    F_2 &= e_5 \wedge (e_6 + \sigma_6 e_1 - \sigma_2 e_2) \wedge e_7
    \wedge e_8~,
  \end{aligned}
\end{equation*}
which moreover are orthogonal.

Finally, if $\lambda_2 = \sigma_2 = 0$, $F$ is given by
\begin{multline*}
  F = e_{1234} + \mu_3 e_{5678}
  + \lambda_5 (e_{2347} + \mu_3 e_{1568})\\
  + \lambda_6 (e_{2348} - \mu_3 e_{1567})
  + \sigma_5 (e_{2345} + \mu_3 e_{1678})
  + \sigma_6 (e_{2346} - \mu_3 e_{1578})~,
\end{multline*}
which can be written as a sum of two orthogonal simple forms
$F= F_1 + \mu_3 F_2$, where
\begin{equation*}
  \begin{aligned}[m]
    F_1 &= (e_1 - \sigma_5 e_5 - \sigma_6 e_6 - \lambda_5 e_7 -
    \lambda_6 e_8) \wedge e_2 \wedge e_3 \wedge e_4\\
    F_2 &= (e_5 + \sigma_5 e_1) \wedge (e_6 + \sigma_6 e_1) \wedge
    (e_7 + \lambda_5 e_1) \wedge (e_8 + \lambda_6 e_1)~.
  \end{aligned}
\end{equation*}

\subsection{Proof for $F \in \Lambda^4\EE^7$}

Choose an orthonormal basis $\{e_1,e_2,\ldots,e_7\}$ for which
$\iota_{12} F = \alpha e_{34} + \beta e_{56} $, where
$\iota_{12}$ means the contraction of $F$ by $e_{12}$.

Suppose that $\alpha$ and $\beta$ are generic.  In this
case, the equation $[\iota_{12} F,F]=0$ says that the only terms in
$F$ which survive are those which are invariant under the maximal
torus of $\SO(5)$, the group of rotations in the five-dimensional
space spanned by $\{e_3,e_4,\ldots,e_7\}$; that is,
\begin{equation*}
  F = \alpha e_{1234} + \beta e_{1256} + \gamma
  e_{3456} ~.
\end{equation*}
Now $[\iota_{23}
F, F] = 0$ implies that $\alpha\beta  = 0$,
violating the condition that $\iota_{12} F$ be generic.

  Non-generic rotations
correspond to (conjugacy classes of) subalgebras of $\fso(5)$ with
rank strictly less than that of $\fso(5)$:
\begin{enumerate}
\item $\fsu(2)$:  $\alpha = \beta \neq 0$; and
\item $\fso(2)$: $\beta = 0$ and $\alpha \neq 0$.
\end{enumerate}
We now go down this list case by case.

\subsubsection{$\fsu(2)$}

In this case $\iota_{12}F = \alpha (e_{34} + e_{56})$.  This means that the most
general solution of $[\iota_{12} F, F]=0$ is given by
\begin{equation*}
  F = \alpha e_{1234} +\alpha  e_{1256} + e_1 \wedge G_1 + e_2 \wedge G_2 + G_3~,
\end{equation*}
where $G_1,G_2 \in \Lambda^3\EE^5$ and $G_3 = \Lambda^4 \EE^5$, where
$\EE^5$ is spanned by $\{e_i\}_{3\leq i\leq 7}$, and where the $G_i$
have zero weight with respect to this $\fsu(2)$ algebra.  A little
group theory shows that $G_1$ and $G_2$ are linear combinations of the
following eight $3$-forms
\begin{equation*}
  e_{347} + e_{567}\qquad e_{347} - e_{567}\qquad e_{357} +
  e_{467}\qquad e_{367} - e_{457}~;
\end{equation*}
 whereas
\begin{equation*}
 G_3 = \mu e_{3456}~.
\end{equation*}
This means that $F$ takes the following form:
\begin{multline*}
  F = \alpha e_{1234} +\alpha e_{1256} + \mu e_{3456}
  + \lambda_1 (e_{1347} + e_{1567})+ \lambda_2 (e_{1347} - e_{1567})
  \\
  +\lambda_3 (e_{1357} +
  e_{1467})+\lambda_4 (e_{1367} - e_{1457})
  + \rho_1 (e_{2347} + e_{2567})
  \\
  + \rho_2 (e_{2347} - e_{2567})
  +\rho_3 (e_{2357} +
  e_{2467})+\rho_4 (e_{2367} - e_{2457})
  ~.
\end{multline*}
Rotating in the anti-selfdual (3456) plane allows us to set $\lambda_3=\lambda_4=0$.  Imposing, for
example, the equation $[\iota_{23}F, F]=0$ and $[\iota_{25}F, F]=0$
tells us that $\lambda_1 =\lambda_2=0$.
This allows us to rotate again in the anti-selfdual (3456) plane to set $\rho_3=\rho_4=0$
and imposing $[\iota_{13}F, F]=0$ and $[\iota_{15}F, F]=0$ to find that $\rho_1 =\rho_2=0$.
The remaining equations imply that $\alpha^2=0$ which is a contradiction.

\subsubsection{$\fso(2)$}

Finally, we consider the case where $\iota_{12} F = \alpha e_{34}$.
The most general $F$ has the form
\begin{equation*}
  F = \alpha e_{1234} + e_1 \wedge G_1 + e_2 \wedge G_2 + G_3~,
\end{equation*}
where $G_1,G_2 \in \Lambda^3\EE^5$ and $G_3 \in \Lambda^4\EE^5$, where
$\EE^5$ is spanned by $\{e_i\}_{3\leq i\leq 7}$.  Such an $F$ will
obey $[\iota_{12} F, F]=0$ if and only if the $G_i$ have zero weights
under the $\fso(2)$ generated by $\iota_{12}F$.  This means that each
of $G_1,G_2$ is a linear combination of the four monomials
\begin{equation*}
  e_{345}\quad e_{346}\quad e_{347}\quad e_{567}\quad
\end{equation*}
Using the freedom to conjugate by the $\SO(3)$ which acts in the
(567) plane, we can write
\begin{equation*}
G_3=\mu e_{3456}~.
\end{equation*}
 So $F$ is
\begin{multline*}
  F = \alpha e_{1234} + \mu e_{3456}
  + \lambda_1 e_{1345} + \lambda_2 e_{1346}
  + \lambda_3 e_{1567}
  \\
+\sigma_1 e_{2345}
  + \sigma_2 e_{2346} + \sigma_3 e_{2567}~.
\end{multline*}
Rotating in the (56)-plane, we can set $\lambda_2=0$. Suppose that
$\mu\not=0$. In this case $[\iota_{36} F, F]=0$ implies that
$\lambda_3=\sigma_3=0$. Next observe that $\iota_{34} F$ is a two-form
in $\EE^4$ spanned by $\{e_1, e_2, e_5, e_6\}$. If $\iota_{34} F$
has rank four then it is the previous case which has led to a contradiction.
If it has rank two, then the statement is shown.

It remains to show the statement for $\mu=0$. In this case, after performing
a rotation in the (56)-plane and setting $\lambda_2=0$, we have
\begin{multline*}
  F = \alpha e_{1234}
  + \lambda_1 e_{1345}
  + \lambda_3 e_{1567}
+\sigma_1 e_{2345}
  + \sigma_2 e_{2346} + \sigma_3 e_{2567}~.
\end{multline*}
One of the $[\iota_{13} F, F]=0$ conditions  implies that $\lambda_1 \sigma_2=0$. If
$\lambda_1=0$, using a rotation in the (56)-plane, we can set $\sigma_2=0$
as well. The conditions $[\iota_{13} F, F]=0$ and $[\iota_{23} F, F]=0$ imply
that $\lambda_3=\sigma_3=0$. Thus
\begin{equation*}
F = \alpha e_{1234}
+\sigma_1 e_{2345}= (\alpha e_1-\sigma_1 e_5) \wedge e_{234}
\end{equation*}
and it is simple. If instead $\sigma_2=0$, using a rotation in the (12)-plane
we can set $\lambda_1=0$. Then an analysis similar to the above yields
that $F$ is simple.

\subsection{Proof for $F \in \Lambda^4\EE^d$  for $d=5,6$}

Choose an orthonormal basis in $\EE^6$ $\{e_1,e_2,\ldots,e_6\}$ for
which $\iota_{12} F = \alpha e_{34} + \beta e_{56} $, where
$\iota_{12}$ means the contraction of $F$ by $e_{12}$.

Suppose that $\alpha$ and $\beta$ are generic.  In this
case, the equation $[\iota_{12} F,F]=0$ says that the only terms in
$F$ which survive are those which are invariant under the maximal
torus of $\SO(4)$, the group of rotations in the five-dimensional
space spanned by $\{e_3,e_4,\ldots,e_6\}$; that is,
\begin{equation*}
  F = \alpha e_{1234} + \beta e_{1256} + \gamma
  e_{3456} ~.
\end{equation*}
Now $[\iota_{23}
F, F] = 0$ implies that $\alpha\beta  = 0$,
violating the condition that $\iota_{12} F$ be generic.

Non-generic rotations correspond to (conjugacy classes of) subalgebras
of $\fso(4)$ with rank strictly less than that of $\fso(4)$:
\begin{enumerate}
\item $\fsu(2)$:  $\alpha = \beta \neq 0$; and
\item $\fso(2)$: $\beta = 0$ and $\alpha \neq 0$.
\end{enumerate}
We now go down this list case by case.

\subsubsection{$\fsu(2)$}

In this case $\iota_{12}F = \alpha (e_{34} + e_{56})$.  This means
that the most general solution of $[\iota_{12} F, F]=0$ is given by
\begin{equation*}
  F = \alpha e_{1234} +\alpha  e_{1256} + e_1 \wedge G_1 + e_2 \wedge
  G_2 + G_3~,
\end{equation*}
where $G_1,G_2 \in \Lambda^3\EE^4$ and $G_3 = \Lambda^4 \EE^4$, where
$\EE^4$ is spanned by $\{e_i\}_{3\leq i\leq 6}$, and where the $G_i$
have zero weight with respect to this $\fsu(2)$ algebra.  A little
group theory shows that $G_1=G_2=0$ and
\begin{equation*}
 G_3 = \mu e_{3456}~.
\end{equation*}
This means that $F$ takes the following form:
\begin{equation*}
  F = \alpha e_{1234} +\alpha e_{1256} + \mu e_{3456} ~.
\end{equation*}
Imposing $[\iota_{23}F, F]=0$ we find that  $\alpha^2=0$ which is a
contradiction.

\subsubsection{$\fso(2)$}

Finally, we consider the case where $\iota_{12} F = \alpha e_{34}$.
The most general $F$ has the form
\begin{equation*}
  F = \alpha e_{1234} + e_1 \wedge G_1 + e_2 \wedge G_2 + G_3~,
\end{equation*}
where $G_1,G_2 \in \Lambda^3\EE^4$ and $G_3 \in \Lambda^4\EE^4$, where
$\EE^4$ is spanned by $\{e_i\}_{3\leq i\leq 6}$.  Such an $F$ will
obey $[\iota_{12} F, F]=0$ if and only if the $G_i$ have zero weights
under the $\fso(2)$ generated by $\iota_{12}F$.  This means that each
of $G_1,G_2$ is a linear combination of the two monomials $e_{345}$
and $e_{346}$, whence
\begin{equation*}
  G_3=\mu e_{3456}~,
\end{equation*}
and
\begin{multline*}
  F = \alpha e_{1234} + \mu e_{3456} + \lambda_1 e_{1345} + \lambda_2
  e_{1346} +\sigma_1 e_{2345} + \sigma_2 e_{2346} ~.
\end{multline*}
Rotating in the (56)-plane, we can set $\lambda_2=0$.  Suppose that
$\mu\not=0$.  Next observe that $\iota_{34} F$ is a two-form in
$\EE^4$ spanned by $\{e_1, e_2, e_5, e_6\}$. If $\iota_{34} F$ has
rank four then it is the previous case which has led to a
contradiction.  If it has rank two, then the statement is shown.

It remains to show the statement for $d=5$.  In this case
\begin{equation*}
  F=\alpha e_{1234}+\beta e_{1534}+\gamma e_{2534}~.
\end{equation*}
The two-form $\iota_{34}F$ has rank two in $\EE^3$ spanned by $\{e_1,
e_2, e_3\}$ and the statement is shown.

\subsection{Proof for $F\in \Lambda^5\EE^{10}$}

We shall not give the details of the proof of the conjecture
in this case. This is because the proof follows closely
that of $F \in \Lambda^5 \EE^{1,9}$ which will be given explicitly
below. The only difference is  certain signs in the
various orthogonality relations that involve the ``time'' direction.
The rest of the proof follows unchanged.

\subsection{Proof for $F \in \Lambda^5 \EE^{1,9}$}

Let us choose a pseudo-orthonormal basis $\{e_0,e_1,\ldots, e_9\}$
with $e_0$ timelike in such a way that the $2$-form $\iota_{012}F$
takes the form
\begin{equation*}
  \iota_{012} F = \alpha e_{34} + \beta e_{56} + \gamma e_{78}.
\end{equation*}
Depending on the values of $\alpha$, $\beta$ and $\gamma$ we have the
same cases as in the case of $d=4$ treated in the previous section.
The most general $F$ can be written as
\begin{multline}
  \label{eq:ansatz5}
  F = \alpha e_{01234} + \beta e_{01256} + \gamma e_{01278} + e_{12}
  \wedge G_0 + e_{02} \wedge G_1 + e_{01} \wedge G_2\\
  + e_0 \wedge H_0 + e_1 \wedge H_1 + e_2 \wedge H_2 + K~,
\end{multline}
where $G_i \in \Lambda^3\EE^7$, $H_i \in \Lambda^4 \EE^7$ and $K \in
\Lambda^5 \EE^7$, where $\EE^7$ is spanned by $\{e_i\}_{3\leq i \leq
9}$.  For all values of $\alpha,\beta,\gamma$, the $2$-form
$\iota_{012} F$ is an element in a fixed  Cartan subalgebra of
$\fso(6)$, and in solving $[\iota_{012}F,F]=0$ we will be determining
which $G_i$, $H_i$ and $K$ have zero weights with respect to this
element.  We will first decompose the relevant exterior powers of
$\EE^7$ in $\fso(6)$ representations.  First of all, notice that
$\EE^7 = \EE^6 \oplus \RR$, where $\EE^6$ is the vector representation
of $\fso(6)$ and $\RR$ is the span of $e_9$.  This means that we can
refine the above decomposition of $F$ and notice that each $G_i$ and
each $H_i$ will be written as follows:
\begin{equation*}
  G_i = L_i + M_i \wedge e_9 \qquad \text{and}\qquad H_i = N_i + P_i
  \wedge e_9~,
\end{equation*}
where $M_i \in\Lambda^2\EE^6$, $L_i,P_i \in\Lambda^3\EE^6$ and
$N_i \in\Lambda^4\EE^6$.  Since $\Lambda^4 \EE^6 \cong
\Lambda^2\EE^6$, we need only decompose $\Lambda^2\EE^6$ and
$\Lambda^3\EE^6$.  Clearly $\Lambda^2\EE^6 \cong \fso(6)$ is nothing
but the $15$-dimensional adjoint representation with three zero
weights corresponding to the Cartan subalgebra, whereas $\Lambda^3
\EE^6$ is a $20$-dimensional irreducible representation having no zero
weights with respect to $\fso(6)$; although of course it many have
zero weights with respect to subalgebras of $\fso(6)$.  Finally, let
us mention that as we saw in the previous section, we will always be
able to choose $K$ to be a linear combination of the monomials
$e_{34569}, e_{34789}, e_{56789}$ by using the freedom to conjugate by
the normaliser of the Cartan subalgebra in which $\iota_{012}F$ lies.

We have different cases to consider depending on the values of
$\alpha$, $\beta$ and $\gamma$ and as in the previous section we can
label them according to the subalgebra of $\fso(6)$ in whose Cartan
subalgebra they lie:
\begin{enumerate}
\item $\fso(6)$: $\alpha$, $\beta$ and $\gamma$ generic;
\item $\fsu(3)$: $\alpha + \beta + \gamma = 0$ but all $\alpha$,
  $\beta$, and $\gamma$ nonzero;
\item $\fsu(2) \times \fu(1)$: $\alpha = \beta \neq \gamma$, but again
  all nonzero;
\item $\fu(1)$ diagonal: $\alpha= \beta= \gamma \neq 0$;
\item $\fso(4)$: $\gamma=0$ and $\alpha \neq \beta$ nonzero;
\item $\fsu(2)$: $\gamma=0$ and $\alpha = \beta \neq 0$; and
\item $\fso(2)$: $\beta = \gamma = 0$ and $\alpha \neq 0$.
\end{enumerate}
We now go down this list case by case.

\subsubsection{$\fso(6)$}

The generic case is easy to discard.  The most general $F$ obeying
$[\iota_{012}F,F]=0$ has 21 free parameters:
\begin{multline*}
  F = \alpha e_{01234} + \beta e_{01256} + \gamma e_{01278} + \mu_1
  e_{34569} + \mu_2 e_{34789} + \mu_3 e_{56789}\\
  + \lambda_1 e_{01349} + \lambda_2 e_{02349} + \lambda_3 e_{12349}
  + \lambda_4 e_{01569} + \lambda_5 e_{02569} + \lambda_6 e_{12569}\\
  + \lambda_7 e_{01789} + \lambda_8 e_{02789} + \lambda_9 e_{12789}
  + \sigma_1 e_{03456} + \sigma_2 e_{03478} + \sigma_3 e_{05678}\\
  + \sigma_4 e_{1345} + \sigma_5 e_{13478} + \sigma_6 e_{15678}
  + \sigma_7 e_{23456} + \sigma_8 e_{23478} + \sigma_9 e_{25678}~.
\end{multline*}
If we now consider the equation $[\iota_{013}F,F]=0$ we see that it
is not satisfied unless either $\alpha$ or $\beta$ are zero, violating
the condition of genericity.

\subsubsection{$\fsu(3)$}

As discussed above, the $\fsu(3)$ zero weights in the
representations $\Lambda^2\EE^6$ and $\Lambda^3\EE^6$ are linear
combinations of the following forms:
\begin{equation*}
  e_{34} \quad e_{56} \quad e_{78} \quad \Omega_1 \quad \Omega_2~,
\end{equation*}
where $\Omega_i$ are defined in equation \eqref{eq:hol3form}.  The
most general $F$ satisfying $[\iota_{012}F,F]=0$ is given by
\begin{multline*}
  F = \alpha (e_{01234} - e_{01278}) + \beta (e_{01256} - e_{01278}) +
  \mu_1 e_{34569} + \mu_2 e_{34789} + \mu_3 e_{56789}\\
  + \lambda_1 e_{01349} + \lambda_2 e_{02349} + \lambda_3 e_{12349}
  + \lambda_4 e_{01569} + \lambda_5 e_{02569} + \lambda_6 e_{12569}\\
  + \lambda_7 e_{01789} + \lambda_8 e_{02789} + \lambda_9 e_{12789}
  + \sigma_1 e_{03456} + \sigma_2 e_{03478} + \sigma_3 e_{05678}\\
  + \sigma_4 e_{1345} + \sigma_5 e_{13478} + \sigma_6 e_{15678}
  + \sigma_7 e_{23456} + \sigma_8 e_{23478} + \sigma_9 e_{25678}\\
  + \rho_1 e_{01} \wedge \Omega_1 + \rho_2 e_{02} \wedge \Omega_1
  + \rho_3 e_{12} \wedge \Omega_1 + \rho_4 e_{01} \wedge \Omega_2
  + \rho_5 e_{02} \wedge \Omega_2 + \rho_6 e_{12} \wedge \Omega_2\\
  - \tau_1 e_{09} \wedge \Omega_1 - \tau_2 e_{19} \wedge \Omega_1
  - \tau_3 e_{29} \wedge \Omega_1 - \tau_4 e_{09} \wedge \Omega_2
  - \tau_5 e_{19} \wedge \Omega_2 - \tau_6 e_{29} \wedge \Omega_2~.
\end{multline*}
There are thus 33 free parameters, which we can reduce to 32 as was
done in the previous section.  Inspection of (some of) the remaining
30239 equations $[\iota_{ijk}F,F]=0$ shows that $\alpha$ and $\beta$
are constrained to obey $\alpha = \pm \beta$, violating the hypothesis
of genericity.

\subsubsection{$\fsu(2) \times \fu(1)$}

We now let $\alpha = \beta$, with the opposite case being related by
an outer automorphism.  As mentioned above $\Lambda^3 \EE^6$ has no
zero weights, whereas those in $\Lambda^2\EE^6$ are linear
combinations of the following forms
\begin{equation*}
  e_{34} + e_{56} \qquad e_{34} - e_{56} \qquad e_{35} + e_{46} \qquad
  e_{36} - e_{45} \qquad e_{78}~.
\end{equation*}
The first and last are the generators of the Cartan subalgebra of
$\fsu(2) \times \fu(1)$ whereas the remaining three are the generators
of the anti-selfdual $\fsu(2) \subset \fso(4)$.  Using the freedom to
conjugate by the anti-selfdual $\fsu(2)$ we will be able to eliminate
two of the free parameters in the expression for $F$, which after this
simplification takes the following form
\begin{multline*}
  F = \alpha (e_{01234} + e_{01256}) + \gamma e_{01278} +
  \mu_1 e_{34569} + \mu_2 e_{34789} + \mu_3 e_{56789}\\
  + \lambda_1 e_{01349} + \lambda_2 e_{02349} + \lambda_3 e_{12349}
  + \lambda_4 e_{01569} + \lambda_5 e_{02569} + \lambda_6 e_{12569}\\
  + \lambda_7 e_{01789} + \lambda_8 e_{02789} + \lambda_9 e_{12789}
  + \sigma_1 e_{03456} + \sigma_2 e_{03478} + \sigma_3 e_{05678}\\
  + \sigma_4 e_{1345} + \sigma_5 e_{13478} + \sigma_6 e_{15678}
  + \sigma_7 e_{23456} + \sigma_8 e_{23478} + \sigma_9 e_{25678}\\
  + \rho_1 (e_{01359} + e_{01469}) + \rho_2 (e_{02359} + e_{02469})
  + \rho_3 (e_{12359} + e_{12469})\\
  + \rho_4 (e_{01369} - e_{01459}) + \rho_5 (e_{02369} - e_{02459}) +
  \rho_6 (e_{12369} - e_{12459})\\
  + \tau_1 (e_{04678} + e_{03578}) + \tau_2 (e_{04578} - e_{03678})
  + \tau_3 (e_{14678} + e_{13578})\\
  + \tau_4 (e_{14578} - e_{13678}) + \tau_5 (e_{24678} + e_{23578}) +
  \tau_6 (e_{24578} - e_{23678})
\end{multline*}
which depends on 33 parameters.  Inspection of the remaining equations
immediately shows that $\alpha\gamma=0$, violating genericity.

\subsubsection{$\fu(1)$ diagonal}

We now let $\alpha = \beta = \gamma$.  As mentioned in the analogous
case in the previous section, $\Lambda^3 \EE^6$ has no zero weights,
whereas those in $\Lambda^2\EE^6$ are linear combinations of the
$\fu(3)$ generators $\omega_i$:
\begin{gather*}
  e_{35} + e_{46} \qquad e_{45} - e_{36} \qquad e_{37} + e_{48} \qquad
  e_{47} - e_{38}\\
  e_{57} + e_{68} \qquad e_{67} - e_{58}\qquad e_{34} \qquad
  e_{56}\qquad e_{78}~.
\end{gather*}
We have the freedom to conjugate by the normaliser of this $\fu(1)$ in
$\fso(6)$, which is precisely $\fu(3)$.  This means that we can
conjugate the $\fu(3)$ generators in the form $K$ in
\eqref{eq:ansatz5} to a Cartan subalgebra of $\fu(3)$.  In summary
the most general $F$ contains 57 parameters and can be written as
\begin{multline*}
  F = \alpha (e_{01234} + e_{01256} + e_{01278}) +
  \mu_1 e_{34569} + \mu_2 e_{34789} + \mu_3 e_{56789}\\
  + \sum_{i=1}^9 \left( \lambda_i e_{01} + \lambda_{9+i} e_{02} +
    \lambda_{18+i} e_{12} \right) \wedge \omega_i\\
  + \sum_{i=1}^9 \left( \sigma_i e_0 + \sigma_{9+i} e_1 +
    \sigma_{18+i} e_2 \right) \wedge \star \omega_i~,
\end{multline*}
where $\star\omega_i \in \Lambda^4\EE^6$ are the Hodge duals of the
$\omega_i$.  Inspection of a few of the remaining equations shows that
they are consistent only if $\alpha = 0$, which violates the
hypothesis.

As in the eight-dimensional case treated in the previous section,
there are no solutions when $\iota_{012} F$ has maximal rank, a fact
which again lacks a simpler explanation.

\subsubsection{$\fso(4)$}

Let $\iota_{012} F = \alpha e_{34} + \beta e_{56}$ with $\alpha$ and
$\beta$ generic.  The condition that $[\iota_{012} F, F] = 0$ means
that $F$ takes the form given by equation \eqref{eq:ansatz5} where
$G_i \in \Lambda^3\EE^7$ are linear combinations of the six monomials
\begin{equation*}
  e_{347}\qquad e_{348}\qquad e_{349}\qquad e_{567}\qquad
  e_{568}\qquad e_{569}~,
\end{equation*}
where the $H_i \in \Lambda^4\EE^7$ are linear combinations of their
duals
\begin{equation*}
  e_{5689}\qquad  e_{5679}\qquad  e_{5678}\qquad  e_{3489}\qquad
  e_{3479}\qquad  e_{3478}~.
\end{equation*}
The $5$-form $K$ is as usual a linear combination of the three
monomials: $e_{34569},e_{34789},e_{56789}$.  In summary, $F$ is given
by the following expression containing 39 free parameters:
\begin{multline*}
  F = \alpha e_{01234} + \beta e_{01256} + \mu_1 e_{34569} + \mu_2
  e_{34789} + \mu_3 e_{56789}\\
  + \lambda_1 e_{01347} + \lambda_2 e_{02347} + \lambda_3 e_{12347}
  + \lambda_4 e_{01348} + \lambda_5 e_{02348} + \lambda_6 e_{12348}\\
  + \lambda_7 e_{01349} + \lambda_8 e_{02349} + \lambda_9 e_{12349}
  + \sigma_1 e_{01567} + \sigma_2 e_{02567} + \sigma_3 e_{12567}\\
  + \sigma_4 e_{01568} + \sigma_5 e_{02568} + \sigma_6 e_{12568}
  + \sigma_7 e_{01569} + \sigma_8 e_{02569} + \sigma_9 e_{12569}\\
  + \rho_1 e_{03478} + \rho_2 e_{13478} + \rho_3 e_{23478}
  + \rho_4 e_{03479} + \rho_5 e_{13479} + \rho_6 e_{23479}\\
  + \rho_7 e_{03489} + \rho_8 e_{13489} + \rho_9 e_{23489}
  + \tau_1 e_{05678} + \tau_2 e_{15678} + \tau_3 e_{25678}\\
  + \tau_4 e_{05679} + \tau_5 e_{15679} + \tau_6 e_{25679}
  + \tau_7 e_{05689} + \tau_8 e_{15689} + \tau_9 e_{25689}~.
\end{multline*}
We can still rotate in the (12) and (78) planes and in this way set
to zero two of the above parameters, say $\sigma_3$ and $\rho_3$,
although we do not gain much from it.  The equations
$[\iota_{ijk}F,F]=0$ have solutions for every $\alpha,\beta$.  Setting
$\alpha=1$ without loss of generality, we find that $\mu_1=0$ and that
all the variables are given in terms of the $\lambda_i$ which remain
unconstrained:
\begin{equation*}
  \begin{aligned}[m]
    \tau_1 &= \mu_3 \lambda_9\\
    \tau_2 &= \mu_3 \lambda_8\\
    \tau_3 &= - \mu_3 \lambda_7\\
    \tau_4 &= - \mu_3 \lambda_6\\
    \tau_5 &= - \mu_3 \lambda_5\\
    \tau_6 &= \mu_3 \lambda_4\\
    \tau_7 &= \mu_3 \lambda_3\\
    \tau_8 &= \mu_3 \lambda_2\\
    \tau_9 &= - \mu_3 \lambda_1\\
  \end{aligned}\qquad
  \begin{aligned}[m]
    \sigma_1 &= - \mu_3 \rho_9\\
    \sigma_2 &= \mu_3 \rho_8\\
    \sigma_3 &= \mu_3 \rho_7\\
    \sigma_4 &= \mu_3 \rho_6\\
    \sigma_5 &= - \mu_3 \rho_5\\
    \sigma_6 &= - \mu_3 \rho_4\\
    \sigma_7 &= - \mu_3 \rho_3\\
    \sigma_8 &= \mu_3 \rho_2\\
    \sigma_9 &= \mu_3 \rho_1\\
  \end{aligned}\qquad
  \begin{aligned}[m]
    \rho_1 &= \lambda_1\lambda_5 - \lambda_2\lambda_4 \\
    \rho_2 &= \lambda_1\lambda_6 - \lambda_3\lambda_4 \\
    \rho_3 &= \lambda_2\lambda_6 - \lambda_3\lambda_5 \\
    \rho_4 &= \lambda_1\lambda_8 - \lambda_2\lambda_7 \\
    \rho_5 &= \lambda_1\lambda_9 - \lambda_3\lambda_7 \\
    \rho_6 &= \lambda_2\lambda_9 - \lambda_3\lambda_8 \\
    \rho_7 &= \lambda_4\lambda_8 - \lambda_5\lambda_7 \\
    \rho_8 &= \lambda_4\lambda_9 - \lambda_6\lambda_7 \\
    \rho_9 &= \lambda_5\lambda_9 - \lambda_6\lambda_8 \\
  \end{aligned}
\end{equation*}
and
\begin{equation*}
  \mu_2 = \lambda_1\lambda_5\lambda_9 - \lambda_3\lambda_5\lambda_7 +
  \lambda_2\lambda_6\lambda_7 + \lambda_3\lambda_4\lambda_8 -
  \lambda_1\lambda_6\lambda_8 - \lambda_2\lambda_4\lambda_9~,
\end{equation*}
subject to one equation
\begin{equation}
  \label{eq:m5so4}
  \beta = \mu_2 \mu_3~.
\end{equation}
Remarkably (perhaps) these equations are precisely the ones that
guarantee that $F$ can be written as a sum of two simple forms
\begin{equation*}
  F = \theta_0 \wedge \theta_1 \wedge \theta_2 \wedge e_3 \wedge e_4 +
  \mu_3 e_5 \wedge e_6 \wedge \theta_7 \wedge \theta_8 \wedge
  \theta_9~,
\end{equation*}
where
\begin{equation*}
  \begin{aligned}[m]
    \theta_0 &= e_0 + \lambda_3 e_7 + \lambda_6 e_8 + \lambda_9 e_9\\
    \theta_1 &= e_1 - \lambda_2 e_7 - \lambda_5 e_8 - \lambda_8 e_9\\
    \theta_2 &= e_2 + \lambda_1 e_7 + \lambda_4 e_8 + \lambda_7 e_9
  \end{aligned}\qquad
  \begin{aligned}[m]
    \theta_7 &= e_7 + \lambda_3 e_0 + \lambda_2 e_1 - \lambda_1 e_2\\
    \theta_8 &= e_8 + \lambda_6 e_0 + \lambda_5 e_1 - \lambda_4 e_2\\
    \theta_9 &= e_9 + \lambda_9 e_0 + \lambda_8 e_1 - \lambda_7 e_2~.
  \end{aligned}
\end{equation*}
Notice moreover that $\theta_i \perp \theta_j$ for $i=0,1,2$ and
$j=7,8,9$, whence the conjecture holds.

\subsubsection{$\fsu(2)$}

Let $\iota_{012} F = \alpha (e_{01234} + e_{01256})$, where we can put
$\alpha = 1$ without loss of generality.  The most general solution of
$[\iota_{012} F, F] = 0$ takes the form \eqref{eq:ansatz5} where $K$
is as usual a linear combination of the three monomials $e_{34569},
e_{34789}, e_{56789}$, the $G_i$ are linear combinations of the
following $3$-forms
\begin{equation*}
  e_{34i} + e_{56i} \qquad e_{34i} - e_{56i} \qquad e_{35i} + e_{46i}
  \qquad e_{36i} - e_{45i} \qquad e_{789}~,
\end{equation*}
where $i=7,8,9$, and the $H_i$ are linear combinations of their
duals.  In total we have 81 free parameters:
\begin{multline*}
  F = e_{01234} + e_{01256} + \mu_1 e_{34569} + \mu_2 e_{34789} +
  \mu_3 e_{56789}\\
  + \lambda_1 e_{01347}
  + \lambda_2 e_{01348}
  + \lambda_3 e_{01349}
  + \lambda_4 e_{01567}
  + \lambda_5 e_{01568}
  + \lambda_6 e_{01569}
  + \lambda_7 e_{01789}\\
  + \lambda_8 (e_{01357} + e_{01467})
  + \lambda_9 (e_{01358} + e_{01468})
  + \lambda_{10} (e_{01359} + e_{01469})\\
  + \lambda_{11} (e_{01367} - e_{01457})
  + \lambda_{12} (e_{01368} - e_{01458})
  + \lambda_{13} (e_{01369} - e_{01459})\\
  + \rho_1 e_{02347}
  + \rho_2 e_{02348}
  + \rho_3 e_{02349}
  + \rho_4 e_{02567}
  + \rho_5 e_{02568}
  + \rho_6 e_{02569}
  + \rho_7 e_{02789}\\
  + \rho_8 (e_{02357} + e_{02467})
  + \rho_9 (e_{02358} + e_{02468})
  + \rho_{10} (e_{02359} + e_{02469})\\
  + \rho_{11} (e_{02367} - e_{02457})
  + \rho_{12} (e_{02368} - e_{02458})
  + \rho_{13} (e_{02369} - e_{02459})\\
  + \sigma_1 e_{12347}
  + \sigma_2 e_{12348}
  + \sigma_3 e_{12349}
  + \sigma_4 e_{12567}
  + \sigma_5 e_{12568}
  + \sigma_6 e_{12569}
  + \sigma_7 e_{12789}\\
  + \sigma_8 (e_{12357} + e_{12467})
  + \sigma_9 (e_{12358} + e_{12468})
  + \sigma_{10} (e_{12359} + e_{12469})\\
  + \sigma_{11} (e_{12367} - e_{12457})
  + \sigma_{12} (e_{12368} - e_{12458})
  + \sigma_{13} (e_{12369} - e_{12459})\\
  + \eta_1 e_{03456}
  + \eta_2 e_{03478}
  + \eta_3 e_{03479}
  + \eta_4 e_{03489}
  + \eta_5 e_{05678}
  + \eta_6 e_{05679}
  + \eta_7 e_{05689}\\
  + \eta_8 (e_{03578} + e_{04678})
  + \eta_9 (e_{03579} + e_{04679})
  + \eta_{10} (e_{03589} + e_{04689})\\
  + \eta_{11} (e_{03678} - e_{04578})
  + \eta_{12} (e_{03679} - e_{04579})
  + \eta_{13} (e_{03689} - e_{04589})\\
  + \phi_1 e_{13456}
  + \phi_2 e_{13478}
  + \phi_3 e_{13479}
  + \phi_4 e_{13489}
  + \phi_5 e_{15678}
  + \phi_6 e_{15679}
  + \phi_7 e_{15689}\\
  + \phi_8 (e_{13578} + e_{14678})
  + \phi_9 (e_{13579} + e_{14679})
  + \phi_{10} (e_{13589} + e_{14689})\\
  + \phi_{11} (e_{13678} - e_{14578})
  + \phi_{12} (e_{13679} - e_{14579})
  + \phi_{13} (e_{13689} - e_{14589})\\
  + \tau_1 e_{23456}
  + \tau_2 e_{23478}
  + \tau_3 e_{23479}
  + \tau_4 e_{23489}
  + \tau_5 e_{25678}
  + \tau_6 e_{25679}
  + \tau_7 e_{25689}\\
  + \tau_8 (e_{23578} + e_{24678})
  + \tau_9 (e_{23579} + e_{24679})
  + \tau_{10} (e_{23589} + e_{24689})\\
  + \tau_{11} (e_{23678} - e_{24578})
  + \tau_{12} (e_{23679} - e_{24579})
  + \tau_{13} (e_{23689} - e_{24589})~.
\end{multline*}

We notice first of all that the equations $[\iota_{ijk}F,F]=0$ imply
that $\lambda_7 = \rho_7 = \sigma_7 = 0$ and after close inspection of
the equations one can see that there are no solutions unless
$\mu_1=0$, which we will assume from now on.

One then must distinguish between two cases, depending on whether or
not $\mu_2$ equals $\mu_3$.  Let us first of all consider the generic
situation $\mu_2 \neq \mu_3$.  One immediately sees that the following
coefficients vanish: $\lambda_i = \rho_i = \sigma_i = \eta_i =
\tau_i = \phi_i = 0$ for $i\geq 8$, leaving $F$ in the following
form
\begin{equation*}
   F = e_{34} \wedge G_1 + e_{56} \wedge G_2~,
\end{equation*}
where
\begin{multline*}
  G_1 = e_{012} + \mu_2 e_{789} + \lambda_1 e_{017} + \lambda_2
  e_{018} + \lambda_3 e_{019} + \rho_1 e_{027} + \rho_2 e_{028} +
  \rho_3 e_{029}\\
  + \sigma_1 e_{127} + \sigma_2 e_{128} + \sigma_3 e_{129} + \eta_2
  e_{078} + \eta_3 e_{079} + \eta_4 e_{089}\\
  + \phi_2 e_{178} + \phi_3 e_{179} + \phi_4 e_{189} + \tau_2 e_{278}
  + \tau_3 e_{279} + \tau_4 e_{289}
\end{multline*}
and
\begin{multline*}
  G_2 = e_{012} + \mu_3 e_{789} + \lambda_4 e_{017} + \lambda_5
  e_{018} + \lambda_6 e_{019} + \rho_4 e_{027} + \rho_5 e_{028} +
  \rho_6 e_{029}\\
  + \sigma_4 e_{127} + \sigma_5 e_{128} + \sigma_6 e_{129} + \eta_5
    e_{078} + \eta_6 e_{079} + \eta_7 e_{089}\\
  + \phi_5 e_{178} + \phi_6 e_{179} + \phi_7 e_{189} + \tau_5 e_{278}
    + \tau_6 e_{279} + \tau_7 e_{289}~.
\end{multline*}
Some of the remaining equations express the $\eta$s, $\phi$s and
$\tau$s in terms of the $\lambda$s, $\rho$s and $\sigma$s:
\begin{equation*}
  \begin{aligned}[m]
    \eta_2 &= \mu_2\sigma_6\\
    \eta_3 &= -\mu_2\sigma_5\\
    \eta_4 &= \mu_2\sigma_4\\
    \eta_5 &= \mu_3\sigma_3\\
    \eta_6 &= -\mu_3\sigma_2\\
    \eta_7 &= \mu_3\sigma_1
  \end{aligned}\qquad
  \begin{aligned}[m]
    \tau_2 &= -\mu_2\lambda_6\\
    \tau_3 &= \mu_2\lambda_5\\
    \tau_4 &= -\mu_2\lambda_4\\
    \tau_5 &= -\mu_3\lambda_3\\
    \tau_6 &= \mu_3\lambda_2\\
    \tau_7 &= -\mu_3\lambda_1
  \end{aligned}\qquad
  \begin{aligned}[m]
    \phi_2 &= \mu_2\rho_6\\
    \phi_3 &= -\mu_2\rho_5\\
    \phi_4 &= \mu_2\rho_4\\
    \phi_5 &= \mu_3\rho_3\\
    \phi_6 &= -\mu_3\rho_2\\
    \phi_7 &= \mu_3\rho_1~;
  \end{aligned}
\end{equation*}
whereas others in turn relate $\lambda_i$, $\rho_i$ and $\sigma_i$ for
$i=4,5,6$  to $\lambda_j$, $\rho_j$ and $\sigma_j$ for $j=1,2,3$:
\begin{footnotesize}
\begin{equation*}
  \begin{aligned}[m]
    \lambda_4 &= \mu_3 (\rho_3\sigma_2 - \rho_2\sigma_3)\\
    \lambda_5 &= \mu_3 (\rho_1\sigma_3 - \rho_3\sigma_1)\\
    \lambda_6 &= \mu_3 (\rho_2\sigma_1 - \rho_1\sigma_2)
  \end{aligned}\qquad
  \begin{aligned}[m]
    \rho_4 &= \mu_3 (\lambda_2\sigma_3 - \lambda_3\sigma_2)\\
    \rho_5 &= \mu_3 (\lambda_3\sigma_1 - \lambda_1\sigma_3)\\
    \rho_6 &= \mu_3 (\lambda_1\sigma_2 - \lambda_2\sigma_1)
  \end{aligned}\qquad
  \begin{aligned}[m]
    \sigma_4 &= \mu_3 (\lambda_2\rho_3 - \lambda_3\rho_2)\\
    \sigma_5 &= \mu_3 (\lambda_3\rho_1 - \lambda_1\rho_3)\\
    \sigma_6 &= \mu_3 (\lambda_1\rho_2 - \lambda_2\rho_1)~.
  \end{aligned}
\end{equation*}
\end{footnotesize}
The remaining independent variables are subject to two final
equations:
\begin{equation}
  \label{eq:m5su2-1}
    \mu_2 = \sum_{\pi\in\fS_3} (-1)^{|\pi|} \lambda_{\pi(1)}
    \rho_{\pi(2)} \sigma_{\pi(3)} \qquad\text{and}\qquad \mu_2 \mu_3 =
    1~,
\end{equation}
where the sum in the first equation is over the permutations of three
letters and weighted by the sign of the permutation.  These equations
guarantee that $G_1$ and $G_2$ are simple forms:
\begin{equation*}
  G_1 = \theta_0 \wedge \theta_1 \wedge \theta_2 \qquad \text{and}
  \qquad
  G_2 = \mu_3 \theta_7 \wedge \theta_8 \wedge \theta_9~,
\end{equation*}
where
\begin{equation*}
  \begin{aligned}[m]
    \theta_0 &= e_0 + \sigma_1 e_7 + \sigma_2 e_8 + \sigma_3 e_9\\
    \theta_1 &= e_1 - \rho_1 e_7 - \rho_2 e_8 - \rho_3 e_9\\
    \theta_2 &= e_2 + \lambda_1 e_7 + \lambda_2 e_8 + \lambda_3 e_9\\[5pt]
    \theta_7 &= e_7 + \sigma_1 e_0 + \rho_1 e_1 - \lambda_1 e_2\\
    \theta_8 &= e_8 + \sigma_2 e_0 + \rho_2 e_1 - \lambda_2 e_2\\
    \theta_9 &= e_9 + \sigma_3 e_0 + \rho_3 e_1 - \lambda_3 e_2~.
  \end{aligned}
\end{equation*}
If we define $\theta_i = e_i$ for $i=3,4,5,6$ then we see that the
$\theta_i$ are mutually orthogonal and hence that
\begin{equation*}
  F = \theta_0 \wedge \theta_1 \wedge \theta_2 \wedge \theta_3 \wedge
  \theta_4 +  \mu_3 \theta_5 \wedge \theta_6 \wedge \theta_7 \wedge
  \theta_8 \wedge \theta_9
\end{equation*}
is a sum of two orthogonal simple forms.

Finally we consider the case $\mu_2 = \mu_3$ which has no solution
unless $\mu_2^2 =1$.  As in the case of four-forms in eight dimensions
treated in the previous section, we will show that we can choose a
frame where the coefficients $\lambda_i$, $\rho_i$ and $\sigma_i$
vanish for $i\geq 8$, thus reducing this case to the generic case
treated immediately above.

Some of the equations $[\iota_{ijk}F,F]=0$ express the $\eta$s,
$\tau$s and $\phi$s in terms of the $\lambda$s, $\rho$s and $\sigma$s,
leaving $F$ in the following form
\begin{scriptsize}
\begin{multline*}
  F = e_{01234} + e_{01256} + \mu_2 (e_{34789} + e_{56789})\\
  + \lambda_1 (e_{01347} - \mu_2 e_{25689}) + \lambda_2 (e_{01348} +
  \mu_2 e_{25679})
  + \lambda_3 (e_{01349} - \mu_2 e_{25678})\\
  + \lambda_4 (e_{01567} - \mu_2 e_{23489}) + \lambda_5 (e_{01568} +
  \mu_2 e_{23479})
  + \lambda_6 (e_{01569} - \mu_2 e_{23478})\\
  + \lambda_8 (e_{01357} + e_{01467} + \mu_2 e_{23589} + \mu_2
  e_{24689})
  + \lambda_9 (e_{01358} + e_{01468} - \mu_2 e_{23579} - \mu_2 e_{24679})\\
  + \lambda_{10} (e_{01359} + e_{01469} + \mu_2 e_{23578} + \mu_2 e_{24678})
  + \lambda_{11} (e_{01367} - e_{01457} + \mu_2 e_{23689} - \mu_2
  e_{24589})\\
  + \lambda_{12} (e_{01368} - e_{01458} - \mu_2 e_{23679} + \mu_2 e_{24579})
  + \lambda_{13} (e_{01369} - e_{01459} + \mu_2 e_{23678} - \mu_2 e_{24578})\\
  + \rho_1 (e_{02347} + \mu_2 e_{15689}) + \rho_2 (e_{02348} - \mu_2
  e_{15679})
  + \rho_3 (e_{02349} + \mu_2 e_{15678})\\
  + \rho_4 (e_{02567} + \mu_2 e_{13489}) + \rho_5 (e_{02568} - \mu_2
  e_{13479})
  + \rho_6 (e_{02569} + \mu_2 e_{13478})\\
  + \rho_8 (e_{02357} + e_{02467} - \mu_2 e_{13589} - \mu_2 e_{14689})
  + \rho_9 (e_{02358} + e_{02468} + \mu_2 e_{13579} + \mu_2 e_{14679})\\
  + \rho_{10} (e_{02359} + e_{02469} - \mu_2 e_{13578} - \mu_2
  e_{14678})
  + \rho_{11} (e_{02367} - e_{02457} - \mu_2 e_{13689} + \mu_2 e_{14589})\\
  + \rho_{12} (e_{02368} - e_{02458} + \mu_2 e_{13679} - \mu_2
  e_{14579})
  + \rho_{13} (e_{02369} - e_{02459} - \mu_2 e_{13678} + \mu_2 e_{14578})\\
  + \sigma_1 (e_{12347} + \mu_2 e_{05689}) + \sigma_2 (e_{12348} -
  \mu_2 e_{05679})
  + \sigma_3 (e_{12349} + \mu_2 e_{05678})\\
  + \sigma_4 (e_{12567} + \mu_2 e_{03489}) + \sigma_5 (e_{12568} -
  \mu_2 e_{03479})
  + \sigma_6 (e_{12569} + \mu_2 e_{03478})\\
  + \sigma_8 (e_{12357} + e_{12467} - \mu_2 e_{03589} - \mu_2
  e_{04689})
  + \sigma_9 (e_{12358} + e_{12468} + \mu_2 e_{03579} + \mu_2 e_{04679})\\
  + \sigma_{10} (e_{12359} + e_{12469} - \mu_2 e_{03578} - \mu_2
  e_{04678})
  + \sigma_{11} (e_{12367} - e_{12457} - \mu_2 e_{03689} + \mu_2 e_{04589})\\
  + \sigma_{12} (e_{12368} - e_{12458} + \mu_2 e_{03679} - \mu_2
  e_{04579}) + \sigma_{13} (e_{12369} - e_{12459} - \mu_2 e_{03678} +
  \mu_2 e_{04578})~.
\end{multline*}
\end{scriptsize}

Let us define the following (anti)selfdual $3$-forms in the (012789)
plane:
\begin{equation*}
  \begin{aligned}[m]
    \omega_0^\pm &= e_{012} \pm \mu_2 e_{789}\\
    \omega_1^\pm &= e_{017} \mp \mu_2 e_{289}\\
    \omega_2^\pm &= e_{018} \pm \mu_2 e_{279}\\
    \omega_3^\pm &= e_{019} \mp \mu_2 e_{278}\\
    \omega_4^\pm &= e_{027} \pm \mu_2 e_{189}
  \end{aligned}\qquad
  \begin{aligned}[m]
    \omega_5^\pm &= e_{028} \mp \mu_2 e_{179}\\
    \omega_6^\pm &= e_{029} \pm \mu_2 e_{178}\\
    \omega_7^\pm &= e_{127} \pm \mu_2 e_{089}\\
    \omega_8^\pm &= e_{128} \mp \mu_2 e_{079}\\
    \omega_9^\pm &= e_{129} \pm \mu_2 e_{078}
  \end{aligned}
\end{equation*}
and the following (anti)selfdual $2$-forms in the (3456) plane:
\begin{equation*}
  \Theta_1^\pm = e_{34} \pm e_{56} \qquad
  \Theta_2^\pm = e_{35} \mp e_{46} \qquad
  \Theta_3^\pm = e_{36} \pm e_{45}~,
\end{equation*}
in terms of which we can rewrite $F$ in a more transparent form:
\begin{equation*}
  F = \Theta_1^+ \wedge \left( \omega_0^+ + \sum_{i=1}^9 \nu_i^+
    \omega_i^+ \right) + \sum_{i=1}^9 \omega_i^- \Psi_i^-~,
\end{equation*}
where the $\Psi_i^-$ are defined by
\begin{equation*}
  \begin{aligned}[m]
    \Psi_1^- &= \nu_1^- \Theta_1^- + \lambda_8 \Theta_2^- +
    \lambda_{11} \Theta_3^-\\
    \Psi_2^- &= \nu_2^- \Theta_1^- + \lambda_9 \Theta_2^- +
    \lambda_{12} \Theta_3^-\\
    \Psi_3^- &= \nu_3^- \Theta_1^- + \lambda_{10} \Theta_2^- +
    \lambda_{13} \Theta_3^-\\
    \Psi_4^- &= \nu_4^- \Theta_1^- + \rho_8 \Theta_2^- +  \rho_{11}
    \Theta_3^-\\
    \Psi_5^- &= \nu_5^- \Theta_1^- + \rho_9 \Theta_2^- +  \rho_{12}
    \Theta_3^-\\
    \Psi_6^- &= \nu_6^- \Theta_1^- + \rho_{10} \Theta_2^- +  \rho_{13}
    \Theta_3^-\\
    \Psi_7^- &= \nu_7^- \Theta_1^- + \sigma_8 \Theta_2^- +
    \sigma_{11} \Theta_3^-\\
    \Psi_8^- &= \nu_8^- \Theta_1^- + \sigma_9 \Theta_2^- +
    \sigma_{12} \Theta_3^-\\
    \Psi_9^- &= \nu_9^- \Theta_1^- + \sigma_{10} \Theta_2^- +
    \sigma_{13} \Theta_3^-~,
  \end{aligned}
\end{equation*}
and where we have introduced the following variables
\begin{equation*}
  \begin{aligned}[m]
    \nu_1^\pm &= \half (\lambda_1 \pm \lambda_4)\\
    \nu_2^\pm &= \half (\lambda_2 \pm \lambda_5)\\
    \nu_3^\pm &= \half (\lambda_3 \pm \lambda_6)
\end{aligned}\qquad
  \begin{aligned}[m]
    \nu_4^\pm &= \half (\rho_1 \pm \rho_4)\\
    \nu_5^\pm &= \half (\rho_2 \pm \rho_5)\\
    \nu_6^\pm &= \half (\rho_3 \pm \rho_6)
\end{aligned}\qquad
  \begin{aligned}[m]
    \nu_7^\pm &= \half (\sigma_1 \pm \sigma_4)\\
    \nu_8^\pm &= \half (\sigma_2 \pm \sigma_5)\\
    \nu_9^\pm &= \half (\sigma_3 \pm \sigma_6)~.
  \end{aligned}
\end{equation*}
Some of the remaining equations $[\iota_{ijk}F,F]=0$ now say that the
nine anti-selfdual $2$-forms $\Psi_i^-$ are collinear.  This means that
by an anti-selfdual rotation in the (3456) plane we can set $\lambda_i
= \rho_i = \sigma_i = 0$ for $i\geq 8$.  We have therefore managed to
reduce this case to the generic case ($\mu_2 \neq \mu_3$) except that
now $\mu_2 = \mu_3$; but this was shown above to verify the
conjecture.

\subsubsection{$\fso(2)$}

Let $\iota_{012} F = \alpha e_{01234}$, where we can put $\alpha = 1$
without loss of generality.  The most general solution of
$[\iota_{012} F, F] = 0$ takes the form \eqref{eq:ansatz5} where the
$K$ is as usual a linear combination of the three monomials
$e_{34569}, e_{34789}, e_{56789}$, the $G_i$ are linear combinations
of the following $3$-forms:
\begin{gather*}
  e_{345}\qquad e_{346}\qquad e_{347}\qquad e_{348}\qquad e_{349}\\
  e_{567}\qquad e_{568}\qquad e_{569}\qquad e_{578}\qquad e_{579}\\
  e_{589}\qquad e_{678}\qquad e_{679}\qquad e_{689}\qquad e_{789}
\end{gather*}
and the $H_i$ are linear combinations of their duals.  The most
general solution to $[\iota_{012}F,F]=0$ has 93 free parameters:
\begin{multline*}
  F = e_{01234} + \mu_{1} e_{34569} + \mu_{2} e_{34789} + \mu_{3}
  e_{56789}\\
  + \lambda_{1} e_{01345} + \lambda_{2} e_{01346} +   \lambda_{3}
    e_{01347} + \lambda_{4} e_{01348} + \lambda_{5} e_{01349}\\
  + \lambda_{6} e_{01567} + \lambda_{7} e_{01568} + \lambda_{8}
  e_{01569} + \lambda_{9} e_{01578} + \lambda_{10} e_{01579}\\
  + \lambda_{11} e_{01589} + \lambda_{12} e_{01678} + \lambda_{13}
  e_{01679} + \lambda_{14} e_{01689} + \lambda_{15} e_{01789}\\
  + \sigma_{1} e_{02345} + \sigma_{2} e_{02346} + \sigma_{3} e_{02347}
  + \sigma_{4} e_{02348} + \sigma_{5} e_{02349}\\
  + \sigma_{6} e_{02567} + \sigma_{7} e_{02568} + \sigma_{8} e_{02569}
  + \sigma_{9} e_{02578} + \sigma_{10} e_{02579}\\
  + \sigma_{11} e_{02589} + \sigma_{12} e_{02678} + \sigma_{13}
  e_{02679} + \sigma_{14} e_{02689} + \sigma_{15} e_{02789}\\
  + \rho_{1} e_{12345} + \rho_{2} e_{12346} + \rho_{3} e_{12347} +
  \rho_{4} e_{12348} + \rho_{5} e_{12349}\\
  + \rho_{6} e_{12567} + \rho_{7} e_{12568} + \rho_{8} e_{12569} +
  \rho_{9} e_{12578} + \rho_{10} e_{12579}\\
  + \rho_{11} e_{12589} + \rho_{12} e_{12678} + \rho_{13} e_{12679} +
  \rho_{14} e_{12689} + \rho_{15} e_{12789}\\
  + \tau_{1} e_{03456} + \tau_{2} e_{03457} + \tau_{3} e_{03458} +
  \tau_{4} e_{03459} + \tau_{5} e_{03467}\\
  + \tau_{6} e_{03468} + \tau_{7} e_{03469} + \tau_{8} e_{03478} +
  \tau_{9} e_{03479} + \tau_{10} e_{03489}\\
  + \tau_{11} e_{05678} + \tau_{12} e_{05679} + \tau_{13} e_{05689} +
  \tau_{14} e_{05789} + \tau_{15} e_{06789}\\
  + \phi_{1} e_{13456} + \phi_{2} e_{13457} + \phi_{3} e_{13458}
  + \phi_{4} e_{13459} + \phi_{5} e_{13467}\\
  + \phi_{6} e_{13468} + \phi_{7} e_{13469} + \phi_{8} e_{13478}
  + \phi_{9} e_{13479} + \phi_{10} e_{13489}\\
  + \phi_{11} e_{15678} + \phi_{12} e_{15679} + \phi_{13}
  e_{15689} + \phi_{14} e_{15789} + \phi_{15} e_{16789}\\
  + \eta_{1} e_{23456} + \eta_{2} e_{23457} + \eta_{3} e_{23458} +
  \eta_{4} e_{23459} + \eta_{5} e_{23467}\\
  + \eta_{6} e_{23468} + \eta_{7} e_{23469} + \eta_{8} e_{23478} +
  \eta_{9} e_{23479} + \eta_{10} e_{23489}\\
  + \eta_{11} e_{25678} + \eta_{12} e_{25679} + \eta_{13} e_{25689} +
  \eta_{14} e_{25789} + \eta_{15} e_{26789}
\end{multline*}

First we consider the case where $\mu_1 \neq 0$.  This means that many
of the parameters must vanish: $\mu_2 = \mu_3 = 0$, $\eta_i = \phi_i
= \tau_i = 0$ for $i\neq 1,4,7$ and  $\lambda_j = \rho_j = \sigma_j =
0$ for $j\neq 1,2,5$.  The resulting $F$ can be written as $F=e_{34}
\wedge G$, where
\begin{multline*}
  G = e_{012} + \mu_{1} e_{569} + \lambda_{1} e_{015} + \lambda_{2}
  e_{016} + \lambda_{5} e_{019}\\
  + \sigma_{1} e_{025} + \sigma_{2} e_{026} + \sigma_{5} e_{029}  +
    \rho_{1} e_{125} + \rho_{2} e_{126} + \rho_{5} e_{129}\\
  + \tau_{1} e_{056} +  \tau_{4} e_{059} + \tau_{7} e_{069}
  + \phi_{1} e_{156} + \phi_{4} e_{159} + \phi_{7} e_{169}\\
  + \eta_{1} e_{256} + \eta_{4} e_{259} + \eta_{7} e_{269}~,
\end{multline*}
where
\begin{equation*}
  \begin{aligned}[m]
    \tau_1 &= \lambda_1\sigma_2-\lambda_2\sigma_1\\
    \tau_4 &= \lambda_1\sigma_5-\lambda_5\sigma_1\\
    \tau_7 &= \lambda_2\sigma_5-\lambda_5\sigma_2
  \end{aligned}\qquad
  \begin{aligned}[m]
    \phi_1 &= \lambda_1\rho_2-\lambda_2\rho_1\\
    \phi_4 &= \lambda_1\rho_5-\lambda_5\rho_1\\
    \phi_7 &= \lambda_2\rho_5-\lambda_5\rho_2
  \end{aligned}\qquad
  \begin{aligned}[m]
    \eta_1 &= \sigma_1\rho_2-\sigma_2\rho_1\\
    \eta_4 &= \sigma_1\rho_5-\sigma_5\rho_1\\
    \eta_7 &= \sigma_2\rho_5-\sigma_5\rho_2~,
  \end{aligned}
\end{equation*}
and subject to the equation
\begin{equation}
  \label{eq:m5u1-1}
  \mu_{1} = \lambda_{5} \rho_{2} \sigma_{1} - \lambda_{2}
  \rho_{5} \sigma_{1} - \lambda_{5} \rho_{1} \sigma_{2} +
  \lambda_{1} \rho_{5} \sigma_{2} + \lambda_{2} \rho_{1}
  \sigma_{5} - \lambda_{1} \rho_{2} \sigma_{5}~,
\end{equation}
which implies that $G$ (and hence $F$) is simple:
\begin{multline*}
    G= (e_0 + \rho_1 e_5 + \rho_2 e_6 + \rho_5 e_9) \wedge (e_1 -
    \sigma_1 e_5 - \sigma_2 e_6 -  \sigma_5 e_9)\\
    \wedge (e_2 + \lambda_1 e_5 + \lambda_2 e_6 + \lambda_5 e_9)~.
\end{multline*}

Let us assume from now on that $\mu_1 = 0$.  If $\mu_2 \neq 0$ then
the same conclusion as above obtains and $F$ is simple.  Details are
the same up to a permutation of the orthonormal basis.  We therefore
assume that $\mu_2=0$.  If $\mu_3=0$ then the following coefficients
vanish: $\eta_i = \phi_i = \tau_i = 0$ for $i\geq 11$ and $\lambda_j
= \rho_j = \sigma_j = 0$ for $j\geq 6$, resulting in $F = e_{34}
\wedge G$, with
\begin{multline*}
  G = e_{012} + \lambda_{1} e_{015} + \lambda_{2} e_{016} +
  \lambda_{3} e_{017} + \lambda_{4} e_{018} + \lambda_{5} e_{019}\\
  + \sigma_{1} e_{025} + \sigma_{2} e_{026} + \sigma_{3} e_{027}
  + \sigma_{4} e_{028} + \sigma_{5} e_{029}\\
  + \rho_{1} e_{125} + \rho_{2} e_{126} + \rho_{3} e_{127} +
  \rho_{4} e_{128} + \rho_{5} e_{129}\\
  + \tau_{1} e_{056} + \tau_{2} e_{057} + \tau_{3} e_{058} +
  \tau_{4} e_{059} + \tau_{5} e_{067}\\
  + \tau_{6} e_{068} + \tau_{7} e_{069} + \tau_{8} e_{078} +
  \tau_{9} e_{079} + \tau_{10} e_{089}\\
  + \phi_{1} e_{156} + \phi_{2} e_{157} + \phi_{3} e_{158}
  + \phi_{4} e_{159} + \phi_{5} e_{167}\\
  + \phi_{6} e_{168} + \phi_{7} e_{169} + \phi_{8} e_{178}
  + \phi_{9} e_{179} + \phi_{10} e_{189}\\
  + \eta_{1} e_{256} + \eta_{2} e_{257} + \eta_{3} e_{258} +
  \eta_{4} e_{259} + \eta_{5} e_{267}\\
  + \eta_{6} e_{268} + \eta_{7} e_{269} + \eta_{8} e_{278} + \eta_{9}
  e_{279} + \eta_{10} e_{289}~,
\end{multline*}
where
\begin{equation}
  \label{eq:linear}
  \begin{aligned}[m]
    \phi_1 &= \lambda_1\rho_2 - \lambda_2\rho_1\\
    \phi_2 &= \lambda_1\rho_3 - \lambda_3\rho_1\\
    \phi_3 &= \lambda_1\rho_4 - \lambda_4\rho_1\\
    \phi_4 &= \lambda_1\rho_5 - \lambda_5\rho_1\\
    \phi_5 &= \lambda_2\rho_3 - \lambda_3\rho_2\\
    \phi_6 &= \lambda_2\rho_4 - \lambda_4\rho_2\\
    \phi_7 &= \lambda_2\rho_5 - \lambda_5\rho_2\\
    \phi_8 &= \lambda_3\rho_4 - \lambda_4\rho_3\\
    \phi_9 &= \lambda_3\rho_5 - \lambda_5\rho_3\\
    \phi_{10} &= \lambda_4\rho_5 - \lambda_5\rho_4
  \end{aligned}\qquad
  \begin{aligned}[m]
    \eta_1 &= \sigma_1\rho_2 - \sigma_2\rho_1\\
    \eta_2 &= \sigma_1\rho_3 - \sigma_3\rho_1\\
    \eta_3 &= \sigma_1\rho_4 - \sigma_4\rho_1\\
    \eta_4 &= \sigma_1\rho_5 - \sigma_5\rho_1\\
    \eta_5 &= \sigma_2\rho_3 - \sigma_3\rho_2\\
    \eta_6 &= \sigma_2\rho_4 - \sigma_4\rho_2\\
    \eta_7 &= \sigma_2\rho_5 - \sigma_5\rho_2\\
    \eta_8 &= \sigma_3\rho_4 - \sigma_4\rho_3\\
    \eta_9 &= \sigma_3\rho_5 - \sigma_5\rho_3\\
    \eta_{10} &= \sigma_4\rho_5 - \sigma_5\rho_4
  \end{aligned}\qquad
  \begin{aligned}[m]
    \tau_1 &= \lambda_1\sigma_2 - \lambda_2\sigma_1\\
    \tau_2 &= \lambda_1\sigma_3 - \lambda_3\sigma_1\\
    \tau_3 &= \lambda_1\sigma_4 - \lambda_4\sigma_1\\
    \tau_4 &= \lambda_1\sigma_5 - \lambda_5\sigma_1\\
    \tau_5 &= \lambda_2\sigma_3 - \lambda_3\sigma_2\\
    \tau_6 &= \lambda_2\sigma_4 - \lambda_4\sigma_2\\
    \tau_7 &= \lambda_2\sigma_5 - \lambda_5\sigma_2\\
    \tau_8 &= \lambda_3\sigma_4 - \lambda_4\sigma_3\\
    \tau_9 &= \lambda_3\sigma_5 - \lambda_5\sigma_3\\
    \tau_{10} &= \lambda_4\sigma_5 - \lambda_5\sigma_4
  \end{aligned}
\end{equation}
subject to the following 10 equations
\begin{equation}
  \label{eq:m5u1-2}
  \sum_{\pi\in\fS_3} (-1)^{|\pi|} \lambda_{\pi(i)} \rho_{\pi(j)}
  \sigma_{\pi(k)} = 0~,
\end{equation}
for $1\leq i < j < k \leq 5$, where the sum is over the permutations
of three letters and weighted by the sign of the permutation.
These equations are precisely the ones which guarantee that
$G$ (and hence $F$) is actually a simple form $G= \theta_0 \wedge
\theta_1 \wedge \theta_2$, with
\begin{equation*}
  \begin{aligned}[m]
    \theta_0 &= e_0 + \rho_1 e_5 + \rho_2 e_6 + \rho_3 e_7 + \rho_4
    e_8 + \rho_5 e_9\\
    \theta_1 &= e_1 - \sigma_1 e_5 - \sigma_2 e_6 - \sigma_3 e_7 -
    \sigma_4 e_8 - \sigma_5 e_9\\
    \theta_2 &= e_2 + \lambda_1 e_5 + \lambda_2 e_6 + \lambda_3 e_7 +
    \lambda_4 e_8 + \lambda_5 e_9~.
  \end{aligned}
\end{equation*}

Finally, if $\mu_3 \neq 0$ all that happens is that we find that the
coefficients which vanish when $\mu_3 = 0$ are given in terms of those
which do not by the following equations:
\begin{equation*}
  \begin{aligned}[m]
    \eta_{15} &= - \mu_3 \lambda_1\\
    \eta_{14} &= \mu_3 \lambda_2\\
    \eta_{13} &= - \mu_3 \lambda_3\\
    \eta_{12} &= \mu_3 \lambda_4\\
    \eta_{11} &= - \mu_3 \lambda_5
  \end{aligned}\qquad
  \begin{aligned}[m]
    \phi_{15} &= \mu_3 \sigma_1\\
    \phi_{14} &= - \mu_3 \sigma_2\\
    \phi_{13} &= \mu_3 \sigma_3\\
    \phi_{12} &= - \mu_3 \sigma_4\\
    \phi_{11} &= \mu_3 \sigma_5
  \end{aligned}\qquad
  \begin{aligned}[m]
    \tau_{15} &= \mu_3 \rho_1\\
    \tau_{14} &= - \mu_3 \rho_2\\
    \tau_{13} &= \mu_3 \rho_3\\
    \tau_{12} &= - \mu_3 \rho_4\\
    \tau_{11} &= \mu_3 \rho_5
  \end{aligned}
\end{equation*}
and
\begin{equation*}
  \begin{aligned}[m]
    \lambda_{15} &= - \mu_3 \eta_1\\
    \lambda_{14} &= \mu_3 \eta_2\\
    \lambda_{13} &= - \mu_3 \eta_3\\
    \lambda_{12} &= \mu_3 \eta_4\\
    \lambda_{11} &= - \mu_3 \eta_5\\
    \lambda_{10} &= \mu_3 \eta_6\\
    \lambda_9 &= - \mu_3 \eta_7\\
    \lambda_8 &= - \mu_3 \eta_8\\
    \lambda_7 &= \mu_3 \eta_9\\
    \lambda_6 &= - \mu_3 \eta_{10}
  \end{aligned}\qquad
  \begin{aligned}[m]
    \rho_{15} &= \mu_3 \tau_1\\
    \rho_{14} &= - \mu_3 \tau_2\\
    \rho_{13} &= \mu_3 \tau_3\\
    \rho_{12} &= - \mu_3 \tau_4\\
    \rho_{11} &= \mu_3 \tau_5\\
    \rho_{10} &= - \mu_3 \tau_6\\
    \rho_9 &= \mu_3 \tau_7\\
    \rho_8 &= \mu_3 \tau_8\\
    \rho_7 &= - \mu_3 \tau_9\\
    \rho_6 &= \mu_3 \tau_{10}
  \end{aligned}\qquad
  \begin{aligned}[m]
    \sigma_{15} &= \mu_3 \phi_1\\
    \sigma_{14} &= - \mu_3 \phi_2\\
    \sigma_{13} &= \mu_3 \phi_3\\
    \sigma_{12} &= - \mu_3 \phi_4\\
    \sigma_{11} &= \mu_3 \phi_5\\
    \sigma_{10} &= - \mu_3 \phi_6\\
    \sigma_9 &= \mu_3 \phi_7\\
    \sigma_8 &= \mu_3 \phi_8\\
    \sigma_7 &= - \mu_3 \phi_9\\
    \sigma_6 &= \mu_3 \phi_{10}~.
  \end{aligned}
\end{equation*}
This implies that $F = F_1 + \mu_3 F_2$, where $F_1$ was shown above
to be simple and $F_2$ is given by
\begin{multline*}
  F_2 = e_{56789} - \eta_{10} e_{01567} + \eta_{9} e_{01568} -
  \eta_{8} e_{01569} - \eta_{7} e_{01578} + \eta_{6} e_{01579}\\
  - \eta_{5} e_{01589} + \eta_{4} e_{01678} - \eta_{3}
  e_{01679} + \eta_{2} e_{01689} - \eta_{1} e_{01789}\\
  + \phi_{10} e_{02567} - \phi_{9} e_{02568} + \phi_{8}
  e_{02569} + \phi_{7} e_{02578} - \phi_{6} e_{02579}\\
  + \phi_{5} e_{02589} - \phi_{4} e_{02678} + \phi_{3}
  e_{02679} - \phi_{2} e_{02689} + \phi_{1} e_{02789}\\
  + \tau_{10} e_{12567} - \tau_{9} e_{12568} + \tau_{8} e_{12569} +
  \tau_{7} e_{12578} - \tau_{6} e_{12579}\\
  + \tau_{5} e_{12589} - \tau_{4} e_{12678} + \tau_{3} e_{12679} -
  \tau_{2} e_{12689} + \tau_{1} e_{12789}\\
  + \rho_{5} e_{05678} - \rho_{4} e_{05679} + \rho_{3} e_{05689} -
  \rho_{2} e_{05789} + \rho_{1} e_{06789}\\
  + \sigma_{5} e_{15678} - \sigma_{4} e_{15679} + \sigma_{3}
  e_{15689} - \sigma_{2} e_{15789} + \sigma_{1} e_{16789}\\
  - \lambda_{5} e_{25678} + \lambda_{4} e_{25679} - \lambda_{3}
  e_{25689} + \lambda_{2} e_{25789} - \lambda_{1} e_{26789}~,
\end{multline*}
where the relations \eqref{eq:linear} hold and the independent
parameters satisfy the same ten equations \eqref{eq:m5u1-2}.  This
then implies that
\begin{equation*}
  F_2 = \theta_5 \wedge \theta_6 \wedge \theta_7 \wedge \theta_8
  \wedge \theta_9~,
\end{equation*}
where
\begin{equation*}
  \begin{aligned}[m]
    \theta_5 &= e_5 + \rho_1 e_0 + \sigma_1 e_1 - \lambda_1 e_2\\
    \theta_6 &= e_6 + \rho_2 e_0 + \sigma_2 e_1 - \lambda_2 e_2\\
    \theta_7 &= e_7 + \rho_3 e_0 + \sigma_3 e_1 - \lambda_3 e_2\\
    \theta_8 &= e_8 + \rho_4 e_0 + \sigma_4 e_1 - \lambda_4 e_2\\
    \theta_9 &= e_9 + \rho_5 e_0 + \sigma_5 e_1 - \lambda_5 e_2~.
  \end{aligned}
\end{equation*}
Finally, we notice that the simple forms $F_1$ and $F_2$ are
orthogonal since so are the one-forms $\theta_i$ (defining $\theta_3 =
e_3$ and $\theta_4=e_4$).  This then concludes the verification of the
conjecture for this case.

\section*{Acknowledgments}

We would like to thank Mohab Abou-Zeid, Robert Bryant, David
Calderbank, Ali Chamseddine, Nigel Hitchin, Blaine Lawson, Felipe
Leitner, Alexandre Pozhidaev, Miles Reid, Dmitriy Rumynin and Wafic
Sabra for useful comments and discussions.

Cette étude fut completée lors d'une visite de JMF à l'IHÉS, qu'il
tient à remercier pour leur invitation et leur soutien ainsi que
d'avoir offert une atmosphère ideale à sa recherche.

JMF is a member of EDGE, Research Training Network HPRN-CT-2000-00101,
supported by The European Human Potential Programme and, in addition,
his research is partially supported by the EPSRC grant GR/R62694/01.

The research of GP is partially supported by the PPARC grants
PPA/G/S/1998/00613 and PPA/G/O/2000/00451 and by the European grant
HPRN-2000-00122.

\bibliographystyle{utphys}
\bibliography{AdS,AdS3,ESYM,Sugra,Geometry,CaliGeo}

\end{document}